\documentclass{article}
\usepackage{hyperref}
\hypersetup{
  colorlinks   = true, 
  urlcolor     = black, 
  linkcolor    = blue, 
  citecolor   = green 
}
\usepackage{geometry}
\usepackage[english]{babel}
\usepackage{amsmath}
\usepackage{graphicx}
\usepackage{ marvosym }
\usepackage[utf8]{inputenc}
\usepackage{mathtools}
\usepackage{geometry}
\usepackage{amsfonts}
\usepackage{amssymb}
\usepackage{amsthm}
\usepackage{thmtools}
\usepackage{t1enc}
\usepackage[titles]{tocloft}
\usepackage{makeidx}
\usepackage{wasysym}
\usepackage{stmaryrd}
\usepackage{calc}  
\usepackage{enumitem} 
\usepackage[refpage]{nomencl}

\usepackage{bookmark}
\usepackage{algpseudocode}
\usepackage{tikz}
\usetikzlibrary{arrows}
\usetikzlibrary{decorations}
\usepackage{float}
\usepackage{mathtools}

\usepackage{enumitem}
\usepackage{linegoal}
\usepackage{calc}

\theoremstyle{plain}
\newtheorem{thm}{Theorem}
\newtheorem{fact}{Fact}

\newtheorem{claim}[thm]{Claim}
\newtheorem{obs}[thm]{Observation}
\newtheorem{prop}[thm]{Proposition}
\newtheorem*{prop*}{Proposition}
\newtheorem{slem}{Sublemma}
\newtheorem{cor}[thm]{Corollary}
\newtheorem{lem}[thm]{Lemma}
\newtheorem{cond}[thm]{Condition}
\newtheorem{conj}[thm]{Conjecture}
\newtheorem*{lem*}{Lemma}
\theoremstyle{definition}
\newtheorem{defn}[thm]{Definition}
\newtheorem*{defn*}{Definition}

\newtheorem{fel*}[thm]{Exercise}

\newtheorem*{megf*}{Observation}
\theoremstyle{remark}

\newtheorem*{rem*}{Remark}

\newenvironment{biz}{\par\noindent{\itshape Proof:}\ }{\rule{1.5ex}{1.5ex}}
\newenvironment{sbiz}{\par\noindent{\itshape Proof:}\ }{\rule{1ex}{1ex}}
\newenvironment{ssbiz}{\par\noindent{\itshape Proof:}\ }{\rule{1ex}{1ex}\ \rule{1ex}{1ex}}

\newenvironment{nbiz}{\par\noindent{\itshape Proof:}\ }{}

\title{Independent and maximal branching packing in infinite matroid-rooted digraphs}
\author{Attila Joó\thanks{MTA-ELTE Egerváry Research Group, Department of Operations Research, Eötvös University, Budapest, Hungary. 
E-mail: {\tt joapaat@cs.elte.hu}.}}
\date{2016}

\begin{document}
\maketitle

\begin{abstract}
 We prove a common generalization of the maximal independent arborescence packing 
 theorem of Cs. Király  
 \cite{kiraly2013maximal}  (which itself is a common generalization of the reachability based arborescence packing result 
 \cite{kamiyama2009arc} and a matroid based arborescence packing result \cite{de2013matroid}) and two of our 
 earlier works about packing branchings in infinite digraphs, namely \cite{joo2015edmonds} and \cite{fenyo2}. 
\end{abstract}

\section{Introduction}
Edmonds' branching packing theorem \cite{edmonds1973edge} has been generalized in several different directions. The most up to date 
survey  in 2016  about these results that we know is in \cite{egres-16-04}. Branching packing problems are mostly investigated in finite 
digraphs but it turned out that in some cases one can relax the finiteness condition of the digraph to some restriction of the 
forward-infinite or backward-infinite directed paths (see \cite{joo2015edmonds} and \cite{fenyo2}). The main result of this paper 
(Theorem \ref{Mf főtétel}) is to give such an infinite generalization of   
\cite{kiraly2013maximal}  which itself is a common generalization of the reachability based arborescence packing result 
\cite{kamiyama2009arc} and a matroid-based arborescence packing result \cite{de2013matroid}. We replace the finiteness of $ D $ 
by some restriction of the behaviour of either its forward-infinite (Condition \ref{Mf végtelen előreutak}) or its backward-infinite 
paths (Condition \ref{Mf végtelen hátrautak}). We also show by examples that 
some obvious further weakenings of our conditions are not possible.

\section{Notations}

We use some basic set theoretic notation. For the power set of $ X $, we write $ \mathcal{P}(X) $. Intersection has higher priority than 
union.  The variables $ \alpha,\beta,\gamma,\xi $ always stand for ordinals. We denote the smallest 
infinite cardinal (i.e. the set of the natural numbers) by $ \omega $.  If $ \kappa $ is a cardinal, then  $ \kappa^{+} $ is its 
successor cardinal.  The restriction of a function 
$ F $ to the subset $ X $ of its domain is denoted by $ F|_{X} $, and $ F[X] $ stands for the image of $ F|_{X} $. We use the 
abbreviation $ 
B-x+y $ for the set $ 
(B\setminus \{ x \})\cup \{ y \} $.
\subsection{Digraphs}
The digraphs $ D=(V,A) $ in this paper may have multiple edges but does not have loops. If $ e\in A $, then  $ \boldsymbol{D-e} $ is an 
abbreviation of $ (V,A\setminus \{ e \}) $. For $ X \subseteq V $, we denote by $ \boldsymbol{D[X]} $ the subdigraph induced by $ X $. If 
the edge $ e $ goes from $ u $ to $ v $, then $\boldsymbol{ \mathsf{tail}(e)}=u $ and $ \boldsymbol{\mathsf{head}(e)}=v $. The set of the 
ingoing and outgoing edges of $ X\subseteq V $ are denoted by $ \boldsymbol{\mathsf{in}_D(X)} $ and $ \boldsymbol{\mathsf{out}_D(X)} $, 
respectively. For a singleton $ 
\{ v \} $, we write $ \mathsf{in}_D(v) $ instead of $ \mathsf{in}_D(\{ v \}) $ and we use this kind of abbreviation in connection with  
singletons in the case of the other 
set-functions as well. 

The paths in this 
paper are directed, repetition of 
vertices is not allowed, and they are finite unless we say explicitly otherwise. We may define paths  
  by the corresponding vertex sequence if parallel edges do not appear there. This sequence determines on ordering $ <_P $ on $ V(P) $.   
  We denote by 
$ \boldsymbol{\mathsf{start}(P) }$ the $ <_P $-smallest  and by $ \boldsymbol{\mathsf{end}(P)} $ the $ <_P $-largest  vertex 
of a path $ P $.  For $ u<_P v $, the subdigraph of $ P $ induced by the elements of the interval $ [u,v] $ is  
denoted by $ \boldsymbol{P[u,v]} $ and called the  segment of $ P $ from $ u $ to $ v $. The initial segments of  
 $ P $ are the segments in the form  $ 
P[\mathsf{start}(P),v] $.  We define terminal segments similarly.   If an initial 
segment of $ P $ is identical to a terminal segment of 
$ Q $, then we may 
join them to a walk and simplify that to a path 
that we call the 
\textbf{concatenation} of $ P $ and $ Q $. We say 
that \textbf{path} $ 
\boldsymbol{P} $ \textbf{goes from} $ \boldsymbol{X} $ \textbf{to} $ \boldsymbol{Y} $ (or shortly $ P $ is a $ 
\boldsymbol{X 
\rightarrow Y} $ path) if  $ \mathsf{start}(P)\in X $ and $ \mathsf{end}(P)\in Y $. Path $ P $ 
goes \textbf{strictly} from $ X $ to $ Y $ 
if exactly the first vertex  of $ P $  is in $ 
X $ and 
exactly the last is in $ Y $ (\textbf{strict} $ X \rightarrow Y $ path). Let $ \boldsymbol{\mathsf{to}_D(X)} $ be the set of those 
vertices from which  $ X $ is reachable by a directed 
path in $ D $.  A path may consist of a single vertex in which case it is a 
\textbf{trivial path}. For a system $ \mathcal{P} $ of paths,  let  $ \boldsymbol{A(\mathcal{P})}=\bigcup_{P\in \mathcal{P}}A(P) $ and we 
denote by  $ 
\boldsymbol{A_{last}(\mathcal{P})} $ the set of the last edges of the (not forward-infinite) paths in $ \mathcal{P} $.

 A digraph $ D $ is called a
\textbf{branching} if it is a directed forest in which every vertex   is 
reachable by a unique path from   $ X:=\{ v\in V(D): \left|\mathsf{in}_D(v)\right|=0 \} $. This $ X $ is the \textbf{root set} of the 
branching. 

\subsection{Infinite matroids}
There were several attempts to extend the notion of matroid by allowing infinite ground sets but keeping the concept of duality. Finally 
in 
\cite{bruhn2013axioms} the authors achieved this goal which made possible the intensive development of the field.    
In this paper we need to use just some very basic facts about infinite matroids. Most of these  are well-known for finite matroids and  
have the same 
proof in the infinite case thus readers with knowledge only about the finite matroids have no disadvantage. In this subsection, we give 
the notations and the facts that we will use in connection with matroids.
    
The pair $ \mathcal{M}=(S,\mathcal{I}) $ is a matroid if $ \mathcal{I}\subseteq \mathcal{P}(S) $ and it satisfies the following axioms.

\begin{enumerate}
\item $ \varnothing\in \mathcal{I} $,
\item $ I\subseteq I'\in\mathcal{I} $ implies $ I\in \mathcal{I} $,
\item  if $ B $ is a $ \subseteq $-maximal element of $ \mathcal{I} $ and  $ I\in \mathcal{I} $ is not maximal, then is an $ i\in 
B\setminus I $ such that $ (I\cup \{ i \})\in \mathcal{I} $,
\item  if $ I\in \mathcal{I} $ and $ I \subseteq X \subseteq S $, then the set $ \{ I'\in \mathcal{I}: I \subseteq I' 
\subseteq X \} $ has a $ \subseteq $-maximal element.
\end{enumerate}
   
The elements of $ \mathcal{I} $ are called \textbf{independent} sets the other subsets of $ S $ are \textbf{dependent}. The $ \subseteq 
$-maximal independent sets (they exist by using axiom 4 with $ I:=\varnothing $ and $ X:=S $) are the \textbf{bases} of the matroid. In 
notation, sometimes we will not distinguish the matroid from its 
ground set unless it would lead misunderstanding. 

\begin{fact}\label{Mf bázis összehas}
If $ B_1 $ and $ B_2 $ are bases of the same matroid and $ B_1\setminus B_2 $ is finite, then $ \left|B_1\setminus B_2\right|= 
\left|B_2\setminus B_1\right| $.
\end{fact}

\noindent It implies that if there is a finite base, then all the bases are finite and have the same size $ \boldsymbol{r(\mathcal{M})} $ 
which is called the \textbf{rank} of $ \mathcal{M} $.   ZFC alone is not able to decide if the bases of a matroid  
have 
necessarily the 
same cardinality. The Generalized Continuum Hypothesis decides the question affirmatively (as shown by D. A. Higgs in 
\cite{higgs1969equicardinality}) but 
it is false under some other set theoretic 
assumptions (proved by N. Bowler and S. Geschke in \cite{bowler2016self}). Hence if there is no finite base, then the rank is simply
 $ \infty $. For $ S' \subseteq S $, 
the 
pair $ 
(S',\mathcal{I}\cap 
\mathcal{P}(S')) $ is a matroid, it is the 
submatroid of $ \mathcal{M} $ that we get by \textbf{restriction} to $ S' $. For $ S'\subseteq S $, we denote by 
$ \boldsymbol{r(S')} $ the rank of the submatroid (corresponding to) $ S' $. 
A $ 
\subseteq 
$-minimal dependent set is called a   \textbf{circuit}.

\begin{fact}\label{Mf minimális dependent van}
A set $ S'\subseteq S $ is dependent if and only if it contains a circuit (which is 
not 
straightforward for an infinite $ S' $ ).
\end{fact}

\begin{fact}
 The relation $ \{ \left\langle x,y  \right\rangle\in S\times S: \exists C\text{ circuit with }x,y\in C \} $ is transitive.
\end{fact}

By adding the diagonals, we may extend the relation above to an equivalence relation. The equivalence classes are called the 
\textbf{components} of the 
matroid. 
 A matroid is called \textbf{finitary} if all of its circuits are finite. In these matroids, an 
infinite set is independent if and only if all of  its finite subsets are independent, in fact this property characterize the finitary 
matroids.

\begin{fact}
Fix a base $ B $ of (the submatroid corresponding to)  $ S' \subseteq S$. Then the subsets $ I $ of $ S\setminus S' $, for which 
$ (I\cup B)\in \mathcal{I} $, forms a matroid on $ S\setminus S' $ and it does not depend on the choice of $ B $. It is the submatroid 
that 
we get by \textbf{contracting} $ S' $.
\end{fact}

\noindent For $ S', S'' \subseteq S $,  we may restrict first $ \mathcal{M} $ to $ S'\cup S'' $ and in the 
resulting matroid contract $ S'' $.  In this case, we denote 
the resulting submatroid by $ \boldsymbol{S'/S''} $. 

If $  \mathcal{M}_\xi=(S_\xi,\mathcal{I}_\xi)\  (\xi<\kappa) $ are matroids with 
pairwise disjoint 
ground sets, then the \textbf{direct sum} $ \mathcal{M}:=\boldsymbol{\bigoplus_{\xi<\kappa}\mathcal{M}_\xi} $ of the matroids $ \{ 
\mathcal{M}_\xi \}_{\xi<\kappa} $ is the matroid on $ \bigcup_{\xi<\kappa}S_\xi $ 
where $ I \in 
\mathcal{I}_{\mathcal{M}} $ if and only if $ (I\cap S_\xi) \in \mathcal{I}_{\xi}$ for all $ \xi<\kappa $. Every matroid is the direct sum 
of its 
components.

 $ i\in S $ is called a
\textbf{loop} if $ \{ i \} $ is a circuit. We denote by $ \boldsymbol{\mathsf{span}(S')} $ the union of $ S' $ and the loops of $ 
S/S'$. 

\begin{fact}
 $ \mathsf{span} $  is a closure operator.
\end{fact}

\begin{fact}[weak circuit elimination]\label{weak circuit elimination}
If $ C_1,C_2 $ are circuits with $ i\in C_1\cap C_2 $, then $ C_1\cap C_2-i $ contains a circuit.
\end{fact}

\begin{cor}
If $ i\in \mathsf{span}(I)\setminus I $  for some independent set $ I $, then there is a unique circuit 
$C\subseteq  I\cup \{ i \} $. Necessarily $ i\in C $ since $ I $ is independent.
\end{cor}
 
\noindent For $ i\in \mathsf{span}(I) $, let us define $ \boldsymbol{C(i,I)}=  
\begin{cases}
\text{the singleton } \{ i \}  & \mbox{if } i\in I,\\
 \text{the unique circuit }  C \text{ above } &\mbox{if }  i\notin I.  
\end{cases}$   
\begin{fact}\label{Mf bázis-becserél}
If $ B $ is a base and $ i\in S\setminus B $, then for any $ j\in C(i,B) $ the set $ B-j+i $ is a base again.
\end{fact}

\begin{cor}\label{Mf bázis-kör}
If $ I $ is independent and $ i\in I\cap\mathsf{span}(J) $  for some $ J\subseteq S $, then there is some $ j\in J $ ($ j=i $ is allowed) 
such 
that $ I-i+j $ is independent.  Furthermore 
$ i $ and $ j $ are in the same component of the matroid and
 if $ I $ was a base, then Fact \ref{Mf bázis összehas} ensures that $ I-i+j 
$ is a base as well.
\end{cor}

One can find a detailed 
survey about the theory of infinite matroids in the Habilitation thesis of N. Bowler 
\cite{nathanhabil}. 

\subsection{Matroid-rooted digraphs}   
We  call a triple $ \boldsymbol{\mathfrak{R}}=(D_{\mathfrak{R}},\mathcal{M}_{\mathfrak{R}},\pi_{\mathfrak{R}}) $ a 
\textbf{matroid-rooted digraph} if $ D_{\mathfrak{R}}=(V,A) $ is a 
digraph, $ 
\mathcal{M}_{\mathfrak{R}}=(S,\mathcal{I})$ is a matroid and  $ \pi_{\mathfrak{R}}: 
S\rightarrow 
\mathcal{P}(V)\setminus \{ \varnothing \} $.   We will omit the subscripts whenever they are clear from the context. For an $ I\in 
\mathcal{I} $ and $ T \subseteq V $ an 
   $ \boldsymbol{(I,T)} $\textbf{-linkage} is a system of edge-disjoint paths $ \{ P_i \}_{i\in I} $ indexed by the elements of $ I $ 
   such 
   that $ P_i $ goes from $ \pi(i) $ to $ T $. In a \textbf{strict} linkage, $ P_i $ goes strictly from $ \pi(i) $ to $ T $.  We say that 
   $ I $ is $ \boldsymbol{T} $\textbf{-linkable} if such a linkage exists. 
 A \textbf{branching packing} 
$ \mathcal{B} $ 
with respect 
to $ \mathfrak{R}  $ is a system of 
edge-disjoint branchings $ \mathcal{B}=\{ \mathcal{B}_i \}_{i\in S} $ in $ D $ where the root set of $ \mathcal{B}_i $ is $ \pi(i) $. A 
branching packing is \textbf{trivial} if none of the branchings in it have any edges. For $ X \subseteq V $, let 
 $ \boldsymbol{\mathcal{S}(X)}=\{ i\in S:\ \pi(i)\cap 
X\neq \varnothing \} $. The 
matroid-rooted digraph is called \textbf{independent} if  $ \mathcal{S}(v)\in \mathcal{I}$ 
for all $ v\in V $. A branching packing is called independent if the matroid-rooted digraph $ 
\mathfrak{B}:=(D,\mathcal{M},\pi_\mathfrak{B}) $ 
is independent 
where $ \pi_{\mathfrak{B}}(i)=V(\mathcal{B}_i) $. Let us denote $ \mathsf{span}\left( \mathcal{S}(\mathsf{to}_D(X)) \right)  $ by   $ 
\boldsymbol{\mathcal{N}(X)} $ (the 
\textbf{need of} $ X $). Clearly $ 
\mathcal{N}(X)=\mathsf{span} \left( \bigcup_{v\in X}\mathcal{N}(v)  \right) $.     
  The branching packing $ \mathcal{B} $ is \textbf{maximal} if for all 
$ v\in V$ the set $ \mathcal{S}_{\mathfrak{B}}(v):= \{ i\in S:\ \pi_\mathfrak{B}(i)\cap 
X\neq \varnothing \} $ spans $ \mathcal{N}(v) $. Hence a branching packing  $ \mathcal{B} $ is 
independent and maximal if and only if $ S_{\mathfrak{B}}(v) $ is a base of $ \mathcal{N}(v) $ for all $ v\in V $. 

\section{Preparations}
\subsection{The linkage condition}
For the existence of a maximal independent branching packing,  the independence of $ \mathfrak{R} $ is obviously necessary since $ 
\mathcal{S}(v)\subseteq \mathcal{S}_{\mathfrak{B}}(v) $  holds for any branching packing $ \mathcal{B} $. The maximality criteria 
leads to  the 
following necessary condition.

\begin{cond}[linkage condition]\label{connected kicsit}
For all  $ v\in V $, there exists a $ (B,v) $-linkage  in $ D $ 
where $ B $ is a base of $ \mathcal{N}(v) $. 
\end{cond}

If we  suppose that $ 
   \mathcal{M} $ and $ D $ are finite, then independence and the linkage condition 
   are enough to ensure 
  the existence of a maximal, independent branching packing as shown by  Cs. Király  in \cite{kiraly2013maximal}. In fact, instead of 
  Condition \ref{connected kicsit} he used the condition ``$ r(\mathcal{S}(X))+\left|\mathsf{in}_D(X)\right|\geq r(\mathcal{N}(X)) $  
  holds for all 
  nonempty 
  $ X \subseteq 
  V $''.  Simple examples show that the literal infinite generalization of this 
  inequality with cardinals  fails to be sufficient in the infinite case. In fact it does not even imply our Condition \ref{connected 
  kicsit} 
  although they are equivalent in the finite case.  

We need a formally stronger (but in fact equivalent) version of Condition \ref{connected kicsit} which 
is more similar with the condition of Cs. Király. 

\begin{cond}\label{connected}
For all nonempty $ X \subseteq V $, there exists a $ (B, X) $-linkage in $ D $  where $ B $ is a base of $ \mathcal{N}(X) $.  
\end{cond}

A linkage above is called a \textbf{linkage for} $ \boldsymbol{X} $ if it is strict and $ B $ contains a base of $ \mathcal{S}(X) $. 
Clearly, 
one can always ensure 
these extra regularity conditions by taking the appropriate segments of the paths and replace some of them with trivial paths.  
Sometimes we will not want to 
deal with these trivial paths. Throwing them away,  the indices of the remaining 
paths form a base of $ \mathcal{N}(X)/\mathcal{S}(X) $. A \textbf{reduced linkage for} $ \boldsymbol{X} $ is a strict $ (B,X) $-linkage 
where $ B $ is a base of $ \mathcal{N}(X)/\mathcal{S}(X) $.

\begin{obs}\label{Mf connect contracted bázisa}
Condition  \ref{connected} is equivalent to demanding the existence of a (reduced) linkage for all nonempty $ X \subseteq V $. 
\end{obs}

\begin{prop}\label{Mf linkage cond equiv}
Condition \ref{connected kicsit} and \ref{connected} are equivalent. 
\end{prop}

\begin{sbiz}
Condition \ref{connected kicsit} is just the restriction of Condition \ref{connected} to the singleton sets $ X=\{ v \}\ (v\in V) $. 
We 
give a 
proof sketch 
for the nontrivial 
direction.  Well-order $ X $ and pick a linkage $ \{ P_i \}_{i\in G_0} $ for the smallest element $ x_0 $ of $ X $.   
Take a 
linkage $ \{ P_i' \}_{i\in B} $ for the following element $ x_1 $. Let $ B'=  \{ i\in B: i\notin \mathsf{span}(G_0) \} $. We claim 
that for $ j\in B' $ the path $ P_j' $ may not have a common edge (or even common vertex) with any path in $ \{ P_i \}_{i\in G_0} $. 
Indeed, if it has, then $ 
j\in \mathcal{N}(x_0) $  and therefore $ j\in 
\mathsf{span}(G_0) $ (since $ G_0 $ spans $ \mathcal{N}(x_0) $), contradicting  the choice $ j\in B' $. But then $ G_1:=G_0\cup 
B' $ spans 
$ \mathcal{N}(\{ x_0,x_1 \}) $  and the path-system $ \{ P_i \}_{i\in G_0}\cup \{ P_i' \}_{i\in B'} $ is edge-disjoint. One can 
finish the 
proof by transfinite 
recursion taking union of the path-systems at limit steps and trim the final system to be independent at the end. 
\end{sbiz}

\subsection{The statement of the main result and feasible extensions}

\noindent We propose the following two possible relaxations of the finiteness of $ D $ and $ \mathcal{M} $.

\begin{cond}\label{Mf végtelen előreutak}
The matroid $ \mathcal{M} $ has finite rank, and for any forward-infinite  path $ P $ the set $ \mathcal{S}(V(P)) $ spans $ 
\mathcal{N}(V(P)) $. 
\end{cond}
\begin{cond}\label{Mf végtelen hátrautak}
 The matroid $ \mathcal{M} $  has at most countably many components, all of which has  finite 
rank. Furthermore, for any backward-infinite  path $ P $ the set $ \mathcal{S}(V(P)) $ spans $ 
\mathcal{N}(V(P)) $. 
\end{cond}

\noindent Now we state our main result.

\begin{thm}\label{Mf főtétel}
If the matroid-rooted digraph  $ \mathfrak{R}=(D,\mathcal{M},\pi) $ satisfies independence, the linkage condition, and either Condition 
\ref{Mf végtelen előreutak} or Condition \ref{Mf végtelen hátrautak}, then there is an independent, maximal branching packing for $ 
\mathfrak{R} $. 
\end{thm}

Instead of dealing with the branchings directly, we introduce the notion of feasible extension of an $ \mathfrak{R} $.  Let $ i_0\in S 
$   
and $ e_0\in A $ such that $ e_0\in \mathsf{out}_D (\pi(i_0)) $ and $  
 \mathcal{S}(\mathsf{head}(e_0))\cup \{ i_0 \}  $ is independent. The matroid-rooted 
digraph  
obtained by $ \boldsymbol{(i_0,e_0)} $\textbf{-extension} from $ \mathfrak{R}=(D,\pi, \mathcal{M}) $ is $ \mathfrak{R}_1:=(D-e_0, 
\mathcal{M}, 
\pi_1) $ where  \[ \pi_1(i)= 
\begin{cases} \pi(i) &\mbox{if } i\neq i_0 \\
\pi(i)\cup \{ \mathsf{head}(e_0) \} & \mbox{if } i=i_0.  
\end{cases}  \]

\noindent This extension is an imitation of giving edge $ e_0 $ to branching $ \mathcal{B}_{i_0} $. A matroid-rooted digraph $ 
\mathfrak{R}' $ is an 
\textbf{extension} of $ \mathfrak{R} $ if there is a transfinite sequence 
(build-sequence) of 
matroid-rooted digraphs  $ \left\langle \mathfrak{R}_\xi: \xi\leq \alpha  
\right\rangle $ (where $ \mathfrak{R}_{\xi}=(D_\xi,\mathcal{M}, \pi_{\xi}) $ and $ D_\xi=(V,A_\xi) $) with the following properties.

\begin{enumerate}
\item $ \mathfrak{R}_0=\mathfrak{R},\ \mathfrak{R}_\alpha=\mathfrak{R}' $,
\item $ \mathfrak{R}_{\beta+1} $ is an $ (i_\beta,e_\beta) $-extension of $ \mathfrak{R}_\beta $ for some $ i_\beta,e_\beta $,
\item for a limit $ \beta $ we have $ \pi_\beta(i)=\bigcup_{\gamma<\beta} \pi_\gamma(i)  $ and $D_\beta:= 
(V,\bigcap_{\gamma<\beta}A_{\gamma} ) $.

\end{enumerate}

\noindent For  an $ \mathfrak{R}' $ extension of $ \mathfrak{R} $, the sequence above is not necessarily unique but the \textbf{order} $ 
\left|\alpha\right| $ of the extension is ($ \left|\alpha\right| 
=\left|A(D_{\mathfrak{R}})\setminus A(D_{\mathfrak{R}'})\right|$).   Let  $ S' \subseteq S $. If $ 
\pi_{\mathfrak{R}}(i)=\pi_{\mathfrak{R}'}(i) $ whenever $ i\notin S' $, then 
we say that $ \mathfrak{R}' $ is an $ \boldsymbol{S'} $\textbf{-extension} of $ \mathfrak{R} $. We define the limit of 
a transfinite sequence of  consecutive extensions in the same way as 
we defined the limit of transfinite sequence of $ (i,e) $-extensions at the limit steps. It is routine to check that for any $ v\in V $ 
and 
any $ 
\mathfrak{R}' $ extension of $ \mathfrak{R} $, we have $ \mathcal{N}_{\mathfrak{R}'}(v)\subseteq \mathcal{N}(v) $.  We call  $ 
\mathfrak{R}' $ a \textbf{feasible extension}  (with respect to 
  $ \mathfrak{R} $) if it satisfies the following condition.

 \begin{cond}\label{Mf fenntartandó csúcsra}
$ \mathfrak{R}' $ is independent and satisfies the linkage condition; furthermore,  $ \mathcal{N}_{\mathfrak{R}'}(v)=\mathcal{N}(v) $ for 
all 
$ 
v\in V $. 
 \end{cond}
 
\noindent In longer terms: $ \mathfrak{R}' $ is independent, and for all $  v\in V $ there is a $ (B,v) $-linkage where $ B $ is a base 
of 
$ 
\mathcal{N}(v) $. 
 It is easy to see that finding a branching packing  $ \{ \mathcal{B}_i \}_{i\in S} $   
       for $ \mathfrak{R} $ is equivalent to finding a feasible extension $ \mathfrak{R}' $ of $ \mathfrak{R} $ such that $ 
       \mathcal{S}_{\mathfrak{R}'}(v) $ is a base of $ \mathcal{N}(v) $ for all $ v\in V $.  Here  $ A(\mathcal{B}_i) $ will consist of 
       those 
       edges $ e $ for which we had an 
      $ (i,e) $-extension in some fixed build-sequence of the extension $ \mathfrak{R}' $.  
    
   Our plan is to construct a build-sequence of such an $ \mathfrak{R}' $ extension. Any extension of an infeasible 
   extension of $ \mathfrak{R} $ is an infeasible extension of $ \mathfrak{R} $, thus every member of the build-sequence needs to be 
   feasible. On the one hand, a feasible extension of a feasible extension of $ \mathfrak{R} $ is clearly a feasible extension of $ 
   \mathfrak{R} $. On the other hand, the limit of feasible extensions is not necessary feasible, therefore it is not enough to  
   to ensure the existence of one single  feasible $ (i,e) $-extension. (In the finite case of course it is enough  since after at most $ 
   \left|A\right| $-many $ (i,e)$-extensions we are done. Furthermore, in this case,  for any independent $ \mathfrak{R} $ 
   that satisfies the 
   linkage condition there exists a feasible $ (i,e) 
   $-extension unless the trivial branching packing is already maximal.)

\subsection{Counterexamples}
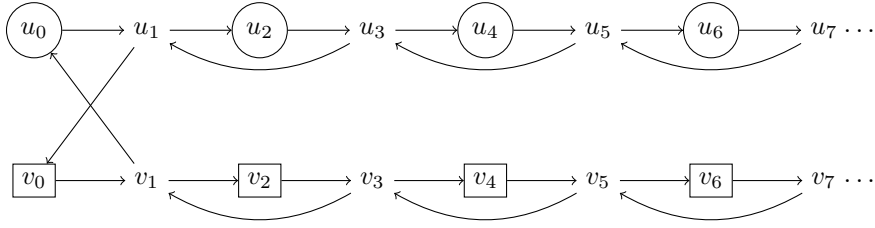
\begin{figure}[H]
\centering

\begin{tikzpicture}

\node (v1) at (-6,1) {$ u_1 $};
\node[rectangle, draw, fill=none] (v2) at (-7.5,-1) {$ v_0 $};
\node (v3) at (-6,-1) {$ v_1 $};
\node[circle, draw, fill=none] (v4) at (-7.5,1) {$ u_0 $};

\draw  (v1) edge[->] (v2);
\draw  (v2) edge[->] (v3);
\draw  (v3) edge[->] (v4);
\draw  (v4) edge[->] (v1);

\node[circle, draw, fill=none] (v5) at (-4.5,1) {$ u_2 $};
\node (v6) at (-3,1) {$ u_3 $};
\node[circle, draw, fill=none] (v7) at (-1.5,1) {$ u_4 $};
\node (v8) at (0,1) {$ u_5 $};
\node[circle, draw, fill=none] (v9) at (1.5,1) {$ u_6 $};
\node (u_1) at (3,1) {$ u_7 $};
\node at (3.5,1) {$ \dots $};

\draw  (v1) edge[->] (v5);
\draw  (v5) edge[->] (v6);
\draw  (v6) edge[->] (v7);
\draw  (v7) edge[->] (v8);
\draw  (v8) edge[->] (v9);
\draw  (v9) edge[->] (u_1);

\node[rectangle, draw, fill=none] (v10) at (-4.5,-1) {$ v_2 $};
\node (v11) at (-3,-1) {$ v_3 $};
\node[rectangle, draw, fill=none] (v12) at (-1.5,-1) {$ v_4 $};
\node (v13) at (0,-1) {$ v_5 $};
\node[rectangle, draw, fill=none] (v14) at (1.5,-1) {$ v_6 $};
\node (u_2) at (3,-1) {$ v_7 $};
\node at (3.5,-1) {$ \dots $};

\draw  (v3) edge[->] (v10);
\draw  (v10) edge[->] (v11);
\draw  (v11) edge[->] (v12);
\draw  (v12) edge[->] (v13);
\draw  (v13) edge[->] (v14);
\draw  (v14) edge[->] (u_2);

\draw  (v11) edge[->,out=-150,in=-30] (v3);
\draw  (v13) edge[->,out=-150,in=-30] (v11);

\draw  (v6) edge[->,out=-150,in=-30] (v1);
\draw  (v8) edge[->,out=-150,in=-30] (v6);

\draw  (u_1) edge[->,out=-150,in=-30] (v8);
\draw  (u_2) edge[->,out=-150,in=-30] (v13);

\end{tikzpicture}
\caption{An independent matroid-rooted digraph  that satisfies the linkage condition  but has no 
feasible $ (i,e) $-extension although $ \mathcal{S}(v) $ is not a base of $ \mathcal{N}(v) $ for any $ v $. $ \mathcal{M} $ is the free 
matroid on $ \{ 0,1 \} $. Elements of $ \pi(0) $ are 
circled and elements of $ \pi(1) $ are in a rectangle in the figure. }\label{Mf kép no feasible extension} 
\end{figure}

As we have already mentioned, independence and the linkage condition are not enough to ensure the existence of an independent maximal 
branching packing. We show this fact by an example (Figure \ref{Mf kép no feasible extension})  where we do not even have a feasible 
$ (i,e) $-extension 
although for any vertex $ v $ the set $ \mathcal{S}(v) $ is not a base of $ \mathcal{N}(v) $.  Let $ V=  
\{ u_n \}_{n<\omega}\cup \{ v_n \}_{n<\omega} $. The edges are $ u_1v_0,\ v_1u_0 $ furthermore, for $ n<\omega $  

 \[  u_nu_{n+1},\ v_nv_{n+1},\  u_{2n+3}u_{2n+1},\  v_{2n+3}v_{2n+1}. \]  Finally take the free matroid on $ \{ 0,1 \} $, 
 let $ \pi(0)=  
\{ u_{2n} \}_{n<\omega} $, and let
$ \pi(1)=  \{ v_{2n} \}_{n<\omega} $. It is routine to check (by using Figure \ref{Mf kép no feasible extension}) that linkage condition 
holds, i.e. every vertex is simultaneously reachable  by edge-disjoint paths from the sets $ \pi(0) $ and $ \pi(1) $. To justify that  
there is no feasible $ (i,e) $-extension, 
 we give for any $ e\in \mathsf{out}_D(\pi(0))$ a vertex set $ X_e $  such that for the 
 $ (0,e) $-extension $ \mathfrak{R}_1 $ we have
 
 \[\mathcal{N}_{\mathfrak{R}_1}(X_e)= \{ 0 \} \subsetneq \{ 0,1 \} = \mathcal{N}(X_e), \]
 which shows the infeasibility. We also do the same for any $ e\in \mathsf{out}_D(\pi(1))$.  For $ n<\omega $, let  $ X_{u_nu_{n+1}}= \{ 
 u_k 
 \}_{n<k<\omega} $ and  let  $ X_{v_nv_{n+1}}= \{ v_k \}_{n<k<\omega} $.
 
 In the example above, $ \pi(0) $ and $ \pi(1) $ are infinite. One can show that if we have a free matroid of arbitrary size and the set 
 $ \pi(i) $ is 
 finite for some $ i $, then there exists an edge $ e $ for which the $ (i,e) $-extension is feasible. Even so, it does not help to 
 construct 
 an independent, maximal branching packing. Indeed, we 
 give an other counterexample with the same matroid where we have $ \pi(0)=\{ u \} $ and $ \pi(1)= \{ v \} $.
  Pick a  2-edge-connected digraph  $ D $ that contains 
  vertices $ u,v $ such that there is no edge-disjoint back and forth paths between $ u $ and $ v $. (Such a digraph exists, even with 
 arbitrary large finite edge-connectivity as we have shown in \cite{joo2016highly}.)  
 From the 
 2-edge-connectivity it follows that every vertex can be reached simultaneously from $ u $ and $ v $ by edge-disjoint paths, thus the 
 linkage condition holds. On the other hand, a  maximal branching packing should contain back and forth paths between $ u $ 
 and $ v $ which do not exist in $ D $.

In the examples above, the structure of the matroid was as simple as possible but the one-way infinite paths do not satisfy any of 
Condition \ref{Mf végtelen előreutak} or Condition \ref{Mf végtelen hátrautak}. Let us give another counterexample (Figure \ref{Mf kép 
finitary 
nem 
elég}) in which, beyond the independence and 
the linkage condition, there is no infinite path at all (not even undirected) and the matroid is just a little bit more complicated than 
what 
Condition \ref{Mf végtelen hátrautak} allows.   

Let $ V= \{ u_n 
\}_{n<\omega}\cup \{ v_n \}_{n<\omega}\cup \{ w \} $ and let $ A=\{ u_nv_n \}_{n<\omega}\cup \{ v_nw \}_{n<\omega}$. The matroid will be 
a countable subset of the vectorspace $ \mathbb{R}^{\omega} $ with the linear independence. We define  
\begin{align*}
 &\mathcal{S}(u_n):=\{(\underbrace{0,\dots,0}_n,1,0,\dots)  
\},\ 
\mathcal{S}(v_n):=\{(\underbrace{0,\dots,0}_n,1,1,0,\dots), (\underbrace{0,\dots,0}_n,1,-1,0,\dots) \},\\   
&\mathcal{S}(w):=\varnothing.
\end{align*} 
  The 
resulting 
matroid-rooted digraph is 
clearly independent. The unique elements of the sets $ \mathcal{S}(u_n) $ form a base of $ \mathcal{N}(w) $, and paths $ u_n, v_n, w \ 
(n<\omega) $ form a linkage for $ w $. 
Considering 
the other vertices,   $ \mathcal{S}(u_n) $ and $ \mathcal{S}(v_n) $ already span $ \mathcal{N}(u_n) $ and $ 
\mathcal{N}(v_n) $ respectively, thus the 
linkage condition 
holds. On the one hand, 
a hypothetical independent and maximal branching packing  may not use any edge of the  from $ u_nv_n $ otherwise it would violate 
independence at $ v_n $. On the other hand we claim that one cannot obtain a base for $ \mathcal{N}(w) $ by 
taking at 
most one element from each $ \mathcal{S}(v_n) $. Indeed, a nontrivial linear combination of such  vectors must have a nonzero component 
other than 
the 0th which ensures that they cannot span $ (1,0,\dots) $. Hence there is no independent and maximal branching packing. 

\begin{figure}[H]
\centering

\begin{tikzpicture}
\node (v1) at (-9.75,0) {$ u_0 $};
\node at (-8.6,0) {$(1,0,\dots)$};
\node (v2) at (-9.75,1.5) {$ v_0 $};
\node at (-8.4,1.2) {$(1,\pm 1,0,\dots)$};

\node (v3) at (-7,0) {$ u_1 $};
\node at (-5.7,0) {$(0,1,0,\dots)$};
\node (v4) at (-7,1.5) {$ v_1 $};
\node at (-5.4,1.3) {$(0,1,\pm 1,0,\dots)$};

\node (v5) at (-3.5,0) {$ u_2 $};
\node at (-2,0) {$(0,0,1,0,\dots)$};
\node (v6) at (-3.5,1.5) {$v_2$};
\node at (-1.6,1.3) {$(0,0,1,\pm 1,0,\dots)$};

\node (v7) at (0.2,0) {$ u_3 $};
\node at (1.9,0) {$(0,0,0,1,0,\dots)$};
\node (v8) at (0.2,1.5) {$ v_3 $};
\node at (2.3,1.3) {$(0,0,0,1,\pm 1,0,\dots)$};

\node (v_1) at (3.5,0.6) {$\dots\dots$};

\node (v9) at (-1.5,3) {$ w $};

\draw  (v1) edge[->] (v2);
\draw  (v3) edge[->] (v4);
\draw  (v5) edge[->] (v6);
\draw  (v7) edge[->] (v8);

\draw  (v2) edge[->] (v9);
\draw  (v4) edge[->] (v9);
\draw  (v6) edge[->] (v9);
\draw  (v8) edge[->] (v9);

\end{tikzpicture}
\caption{An independent matroid-rooted digraph $ \mathfrak{R}=(D,\mathcal{M},\pi) $ that satisfies the linkage condition. Furthermore $ D 
$ 
does 
not contain even undirected infinite paths and $ \mathcal{M} $ 
is countable and finitary   but there is  no 
independent, 
maximal branching packing. For every vertex $ v $, we listed the elements of $ \mathcal{S}(v) $ next to $ v $.} \label{Mf kép finitary 
nem 
elég}
\end{figure}
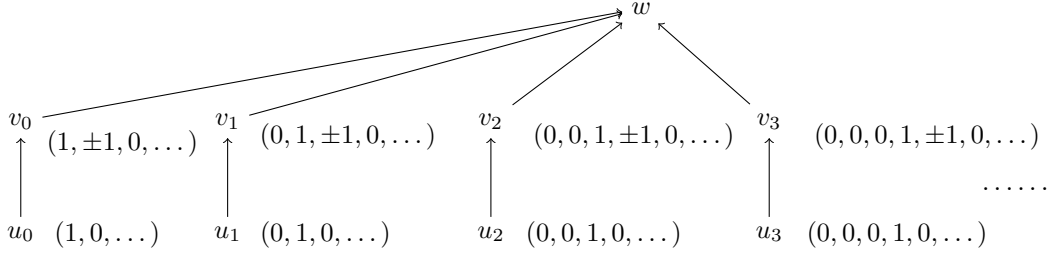

There is an asymmetry in the matroid restriction part of Condition \ref{Mf végtelen előreutak} and Condition \ref{Mf végtelen 
hátrautak}. In our last example, we show that one cannot replace the ``$ \mathcal{M} $ have  finite rank '' part of Condition \ref{Mf 
végtelen előreutak} by the condition that $ \mathcal{M} $ has countably many components all of which has a finite rank. Let be $ V=\{ t 
\}\cup \{ (m,n)\in \omega\times \omega: m \leq n \} $. The  set $ A $ consists of the following edges  (see Figure \ref{Mf ellenpélda 
forwardnál gáz van}). For all $ m,n<\omega $, for which it makes sense

\begin{enumerate}
\item infinitely many parallel edges from $ (m,n+1) $ to $ (m,n) $,
\item  edge from $ (m,n) $ to $ (m+1,n) $,
\item edge from $ (2m+2,n) $ to $ (2m,n) $,
\item edge from $ (m,m) $ to $ t $,
\item edge from $ t $ to $ (2m+1,n) $ (\emph{not in the figure!}).
\end{enumerate}

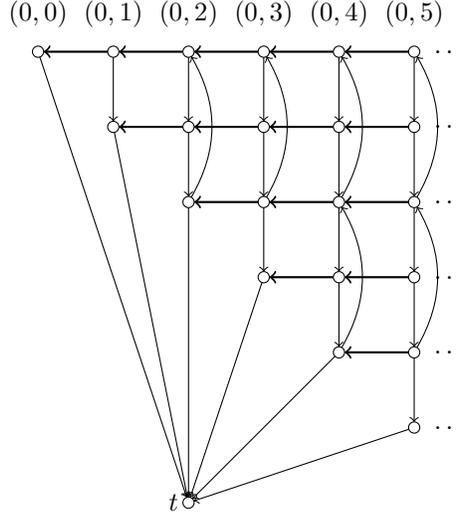
\begin{figure}[H]
\centering

\begin{tikzpicture}
\node[circle, fill=white, draw, outer sep=0pt, inner sep=1.5pt ] (v1) at (-4.5,1.5) {};
\node[circle, fill=white, draw, outer sep=0pt, inner sep=1.5pt ] (v8) at (-3.5,1.5) {};
\node[circle, fill=white, draw, outer sep=0pt, inner sep=1.5pt ] (v9) at (-2.5,1.5) {};
\node[circle, fill=white, draw, outer sep=0pt, inner sep=1.5pt ] (v11) at (-1.5,1.5) {};
\node[circle, fill=white, draw, outer sep=0pt, inner sep=1.5pt ] (v14) at (-0.5,1.5) {};
\node[circle, fill=white, draw, outer sep=0pt, inner sep=1.5pt ] (v18) at (0.5,1.5) {};
\node[circle, fill=white, draw, outer sep=0pt, inner sep=1.5pt ] (v3) at (-3.5,0.5) {};
\node[circle, fill=white, draw, outer sep=0pt, inner sep=1.5pt ] (v10) at (-2.5,0.5) {};
\node[circle, fill=white, draw, outer sep=0pt, inner sep=1.5pt ] (v12) at (-1.5,0.5) {};
\node[circle, fill=white, draw, outer sep=0pt, inner sep=1.5pt ] (v15) at (-0.5,0.5) {};
\node[circle, fill=white, draw, outer sep=0pt, inner sep=1.5pt ] (v19) at (0.5,0.5) {};
\node[circle, fill=white, draw, outer sep=0pt, inner sep=1.5pt ] (v4) at (-2.5,-0.5) {};
\node[circle, fill=white, draw, outer sep=0pt, inner sep=1.5pt ] (v13) at (-1.5,-0.5) {};
\node[circle, fill=white, draw, outer sep=0pt, inner sep=1.5pt ] (v16) at (-0.5,-0.5) {};
\node[circle, fill=white, draw, outer sep=0pt, inner sep=1.5pt ] (v20) at (0.5,-0.5) {};
\node[circle, fill=white, draw, outer sep=0pt, inner sep=1.5pt ] (v5) at (-1.5,-1.5) {};
\node[circle, fill=white, draw, outer sep=0pt, inner sep=1.5pt ] (v17) at (-0.5,-1.5) {};
\node[circle, fill=white, draw, outer sep=0pt, inner sep=1.5pt ] (v21) at (0.5,-1.5) {};
\node[circle, fill=white, draw, outer sep=0pt, inner sep=1.5pt ] (v6) at (-0.5,-2.5) {};
\node[circle, fill=white, draw, outer sep=0pt, inner sep=1.5pt ] (v22) at (0.5,-2.5) {};
\node[circle, fill=white, draw, outer sep=0pt, inner sep=1.5pt ] (v7) at (0.5,-3.5) {};
\node[circle, fill=white, draw, outer sep=0pt, inner sep=1.5pt ] (v2) at (-2.5,-4.5) {};
\draw  (v1) edge[->] (v2);
\draw  (v3) edge[->] (v2);
\draw  (v4) edge[->] (v2);
\draw  (v5) edge[->] (v2);
\draw  (v6) edge[->] (v2);
\draw  (v7) edge[->] (v2);
\draw  (v8) edge[->] (v3);
\draw  (v9) edge[->] (v10);
\draw  (v10) edge[->] (v4);
\draw  (v11) edge[->] (v12);
\draw  (v12) edge[->] (v13);
\draw  (v13) edge[->] (v5);
\draw  (v14) edge[->] (v15);
\draw  (v15) edge[->] (v16);
\draw  (v16) edge[->] (v17);
\draw  (v17) edge[->] (v6);
\draw  (v18) edge[->] (v19);
\draw  (v19) edge[->] (v20);
\draw  (v20) edge[->] (v21);
\draw  (v21) edge[->] (v22);
\draw  (v22) edge[->] (v7);

\node at (-2.7,-4.5) {$t$};

\draw  (v18) edge[thick, ->] (v14);
\draw  (v14) edge[thick, ->] (v11);
\draw  (v11) edge[thick, ->] (v9);
\draw  (v9) edge[thick, ->] (v8);
\draw  (v8) edge[thick, ->] (v1);
\draw  (v19) edge[thick, ->] (v15);
\draw  (v15) edge[thick, ->] (v12);
\draw  (v12) edge[thick, ->] (v10);
\draw  (v10) edge[thick, ->] (v3);
\draw  (v20) edge[thick, ->] (v16);
\draw  (v16) edge[thick, ->] (v13);
\draw  (v13) edge[thick, ->] (v4);
\draw  (v21) edge[thick, ->] (v17);
\draw  (v17) edge[thick, ->] (v5);
\draw  (v22) edge[thick, ->] (v6);
\node at (-4.5,2) {$(0,0)$};
\node at (-3.5,2) {$(0,1)$};
\node at (-2.5,2) {$(0,2)$};
\node at (-1.5,2) {$(0,3)$};
\node at (-0.5,2) {$(0,4)$};
\node at (0.5,2) {$(0,5)$};

\node at (1,1.5) {$\dots$};
\node at (1,0.5){$\dots$};
\node at (1,-0.5) {$\dots$};
\node at (1,-1.5) {$\dots$};
\node at (1,-2.5) {$\dots$};
\node at (1,-3.5) {$\dots$};

\draw  (v4) edge[->,out=60,in=-60] (v9);
\draw  (v13) edge[->,out=60,in=-60] (v11);
\draw  (v16) edge[->,out=60,in=-60] (v14);
\draw  (v20) edge[->,out=60,in=-60] (v18);
\draw  (v6) edge[->,out=60,in=-60] (v16);
\draw  (v22) edge[->,out=60,in=-60] (v20);
\end{tikzpicture}
\caption{An illustration  that in Condition \ref{Mf végtelen előreutak} one cannot replace the restriction of the matroid by the 
weaker  restriction of Condition \ref{Mf végtelen hátrautak}. The outgoing edges of $ t $ (a single edge to each vertex in an odd row) 
are 
not on the 
figure. The thick horizontal edges stand for infinitely many parallel edges.}\label{Mf ellenpélda 
forwardnál gáz van}
\end{figure}

 Observe that after the deletion of $ t $ just finitely many vertices are reachable from any vertex, which shows that there is no 
 forward-infinite path in $ D:=(V,A) $. Let $ \mathcal{M} $ be the free matroid on $ \omega $ and let $ \pi(n)=\{ (0,n) \} $. It is easy 
 to 
 check (using Figure \ref{Mf ellenpélda forwardnál gáz van}) that $ \mathcal{N}(v)=\omega $ for all $ v\in V $ and the linkage condition 
 holds. We have to show that there are 
 no edge-disjoint 
 spanning branchings  with the prescribed root sets. Suppose to the contrary that there is and fix one, say $ \mathcal{B}=\{ 
 \mathcal{B}_n \}_{n<\omega} $.  The only possibility for $ \mathcal{B}_0 $ to reach $ t $ is to use the edge $ ((0,0),t) $. Suppose that 
 we already know for some $ 0<N $ that $ \mathcal{B}_n $ contains the path $ (0,n),(1,n),\dots,(n,n),t $ whenever $ n< N $. By using just 
 the remaining 
 edges, $ t $ is no longer reachable from columns  $0,\dots,N-1 $. It easy to check (using Figure \ref{Mf ellenpélda forwardnál gáz 
 van})  
 that for $ \mathcal{B}_{N} $  the path   
 $ (0,N),(1,N),\dots,(N,N),t $ is the only possible option to reach $ t $. On the other hand after the deletion of the edges of these 
 paths 
 for all 
 $ n $ 
 the 
 vertices $ \{ 
 (0,n): 1\leq n<\omega \} $ are no longer reachable from $ \{ (0,0),t \} $. This prevents $ \mathcal{B}_0 $ from being a spanning 
 branching 
 rooted at $ (0,0) $ which is a contradiction. 
 
 \section{Duality and the characterisation of the infeasible $ \boldsymbol{(i,e)} $-extensions} 
 
Assume that the linkage condition and independence hold for $ \mathfrak{R} $ and let us focus first just on a single  $ (i_0,e_0) 
$-extension 
$ 
\mathfrak{R}_1 
$ of $ \mathfrak{R} $. We cannot ruin 
 independence in this extension, as it is built into the 
    definition of the $ (i,e) $-extension.  If  for some nonempty $ X\subseteq V $ any linkage for $ X $ necessarily  uses all the 
    ingoing edges of $ X $, then we 
     call $ X $ \textbf{tight} (with respect to $ \mathfrak{R} $). 
     If 
     $ X $ is tight and $ 
     i_0\in\mathsf{span}(\mathcal{S}(X)) $, then   $ X $ is called $ \boldsymbol{i_0}$\textbf{-dangerous}. We claim  that if $ e_0 $ 
 is an ingoing edge of 
 an $ i_0 $-dangerous set $ X $, then the  $ (i_0,e_0) $-extension is infeasible. On the one hand, $ 
 i_0\in \mathsf{span}(\mathcal{S}(X)) $ 
 implies that $ \mathsf{span}(\mathcal{S}_{\mathfrak{R}_1}(X)))=\mathsf{span}(\mathcal{S}(X)) $ and hence   $ 
 \mathcal{N}(X)/\mathcal{S}_{\mathfrak{R}_1}(X)=\mathcal{N}(X)/\mathcal{S}(X) $. On the other hand, by the tightness of $ X $ (with 
 respect to $ 
 \mathfrak{R} $)  any $ (B,X) $-linkage   where $ B $ is a base of $ 
 \mathcal{N}(X)/\mathcal{S}(X)= \mathcal{N}(X)/\mathcal{S}_{\mathfrak{R}_1}(X) $  
 uses all the ingoing edges of $ X $ including $ e_0 $  thus there is no more 
 a desired  linkage for $ X $ with respect to  $ \mathfrak{R}_1 $.  It will turn out that surprisingly this is the only possible 
 reason for the infeasibility of an $ (i_0,e_0) $-extension.

 \noindent In the finite case, one can justify this easily in the following way. We use without proof that if $ 
 \mathcal{M} $ has finite rank, then the consequence 
 
  \[\mathfrak{R}'\text{ is independent and } r(\mathcal{S}_{\mathfrak{R}'}(X))+\left|\mathsf{in}_{D_{\mathfrak{R}'}}(X)\right|\geq 
  r(\mathcal{N}(X))\text{ for all nonempty } X 
  \subseteq V\]  
 of Condition \ref{Mf fenntartandó csúcsra}  is actually  equivalent with it. 
 Furthermore, tightness of $ X $ is equivalent with the fact  equality holds for $ X $ in the inequality above. (Of course in the 
 finite case we do not need to know this equivalence or anything about our Condition \ref{Mf fenntartandó csúcsra} at all. One can simply 
 define tightness 
 based on the inequality.)
 
  If the   $  (i_0,e_0)  $-extension 
 is infeasible in the  
 finite case and $ X^{*} $ is a violating set with respect to the resulting $ \mathfrak{R}_1 $, then the extension necessarily reduces 
 the number of 
 ingoing edges of $ X^{*} $ (i.e. $ \left|\mathsf{in}_D(X^{*})\right|=\left|\mathsf{in}_{D-e_0}(X^{*})\right|+1 $ and hence $ e_0\in 
 \mathsf{in}_D(X^{*}) $) but does not increase the rank of 
 the submatroid corresponding to $ X^{*} $ 
 (i.e. $ r(\mathcal{S}_{\mathfrak{R}_1}(X^{*}))=r(\mathcal{S}(X^{*})) $)  thus $i_0\in \mathsf{span}(\mathcal{S}(X^{*})) $, furthermore 
 there must be 
 originally equality for $ X^{*} 
 $. Summarizing these  we obtain that  $ e_0 $ is an ingoing edge of the $ i_0 $-dangerous set $ X^{*} $. 
 
 As we mentioned, the same 
 characterisation of infeasible extensions 
 remains true in the general case, although we need to use more complex arguments to prove it. The rest of the section contains this 
 proof and the corresponding preparations.\\ 
 
  A set $ X\subseteq V $ is called 
   $ \boldsymbol{t} $\textbf{-good} for some $ t\in V $ if there is a system of edge-disjoint paths $ \{ P_b \}_{b\in B}\cup \{ P_e 
   \}_{e\in 
   \mathsf{in}_D(X)} 
   $ in $ D[X] $ such that $ B $ is a base of $ \mathcal{S}(X) $ and  $ \{ P_b \}_{b\in B}  $ is a $ (B,t) $-linkage  and $ P_e $ goes 
   from 
   $ \mathsf{head}(e) $ to $ t $.

  \begin{defn}[complementarity conditions]\label{complementarity conditions}
   The \textbf{complementarity conditions} for an $ (I,t) $-linkage $ \{ P_i \}_{i\in I} $ and a vertex set $ 
  X\ni t $ are the 
  following.
  
  \begin{enumerate}
  \item $ I_{in}:=I\cap \mathcal{S}(X) $ is a base of $ \mathcal{S}(X) $,
  \item  paths $ \{ P_i \}_{i\in I_{in}} $ lie in $ D[X] $,
  \item for $ i\in I\setminus \mathcal{S}(X)=:I_{out} $ we have $ \left|A(P_i)\cap \mathsf{in}_D(X)\right|=1 $,
  \item $ \bigcup_{i\in I_{out}}A(P_i) \supseteq \mathsf{in}_D(X)$.
  \end{enumerate} 
  \end{defn}
\noindent For $ i\in I_{out}$, let us denote by $ e_i $  the first edge of $ P_i $ 
  that enters  $ 
  X $.  Note that if the complementarity conditions hold for $ \mathcal{P} $ and $ X $, then $ X $ is $ t $-good, as shown  by the 
  paths   \[ 
  \{ P_i 
  \}_{i\in I_{in}}\cup \{  
P_i[\mathsf{head}(e_i),t] \}_{i\in 
  I_{out}}. \]  One can  
  replace the 
  conditions 2,3,4 by the single  condition   \[ A_{last} (\{ P_i[\mathsf{start}(P_i),\mathsf{head}(e_i)] \}_{i\in 
  I_{out}})=\mathsf{in}_D(X).  \]

\begin{conj}
 We always have some  $ (I,t) 
$-linkage $ \mathcal{P} $ and an  $ X\ni t $ such that $ \mathcal{P} $ and $ X $ satisfy the 
complementarity conditions.
 \end{conj}

Note that for the free matroid  this conjecture is a reformulation of the famous Infinite Menger theorem \cite{aharoni2009menger} of 
Aharoni and Berger. On the other hand, in the finite case much more general versions are true  (see for example 
\cite{lawler1982computing}).    
    
  \begin{claim}\label{Mf largest t-good}
  There exists a  $ \subseteq $-largest $ t $-good set. 
  \end{claim}
  \begin{ssbiz}
  First of all we always have a smallest $ t $-good set, namely $ \{ t \} $.

  \begin{prop}\label{Mf t-good lánc unió}
  For any $ \subseteq $-increasing, nonempty chain $ \left\langle X_\beta: \beta<\alpha  \right\rangle $  of $ t $-good sets, $ 
    \bigcup_{\beta<\alpha} X_\beta $ is $ t $-good.
  \end{prop}
  
   \begin{sbiz}
   Note that the definition of $ t $-goodness is equivalent if we 
  demand a  generator system $ G $ (a set that contains a base)  instead of a base $ B $ of $ \mathcal{S}(X) $.  We define for all $ 
  \beta 
  \leq \alpha $ a path-system $ 
  \mathcal{P}_\beta $ that shows the $ t $-goodness of $ X_\beta $. Let $ \mathcal{P}_0 $ be an arbitrary system that witnesses the $ t 
  $-goodness of $ X_0 $. If some  $ \mathcal{P}_\beta=\{ P_g 
  \}_{g\in G_\beta}\cup \{P_e: e\in 
  \mathsf{in}_D(X_{\beta}) \} $ has been defined, then we obtain $ 
  \mathcal{P}_{\beta+1} $ in the following way. Let $ \mathcal{P}=\{ P_g' \}_{g\in G'}\cup \{P_e': e\in \mathsf{in}_D(X_{\beta+1}) \} $ 
  be 
  an arbitrary linkage that shows the $ t $-goodness of $ X_{\beta+1} $. Throw away the elements of $ G' $ that are spanned by $ 
  G_\beta $ and take the union of $ G_\beta $ and the reminder of $ G' $ to obtain $ G_{\beta+1} $. For $ g\in G_\beta $, we keep the 
  path 
  $ P_g $ unchanged. Observe that for $ g\in G_{\beta+1}\setminus G_\beta $, the path $ P_g'\in \mathcal{P} $ may not start inside $ 
  X_\beta 
  $, because then $ G_\beta $ would span $ g $ since $ G_\beta $ is a generator for $ \mathcal{S}(X_\beta) $. For $ g\in 
  G_{\beta+1}\setminus 
  G_\beta $, let $ e_g $ be the first   edge 
  of $ P_g' $ that enters $ X_\beta $. We obtain $ P_g $ as a concatenation of paths $ P_g'[\mathsf{start}(P_g'), 
  \mathsf{head}(e_{g})] 
  $ and $ 
  P_{e_g}\in \mathcal{P}_{\beta} $.  We do the same terminal segment replacement process with all the  paths $ \{P_e': e\in 
  \mathsf{in}_D(X_{\beta+1}) \} $ as 
  well for 
  getting $ \{P_e: e\in \mathsf{in}_D(X_{\beta+1}) \} $. Note that the resulting system $ \mathcal{P}_{\beta+1} $ is really 
  edge-disjoint.

  Let $ \beta\leq \alpha $ be a limit ordinal. Observe that $ G_\beta:=\bigcup _{\gamma<\beta}G_\gamma $ is a generator system 
  for $ \mathcal{S}(X_\beta) $. Indeed, if $ i\in \mathcal{S}(X_\beta) $, then $ i\in \mathcal{S}(X_\gamma) $ for some $ \gamma<\beta $ 
  hence $i\in 
  \mathsf{span}(G_\gamma)\subseteq \mathsf{span}(G_\beta) $. If $ e\in \mathsf{in}_D(X_\beta) $, then $ e\in \mathsf{in}_D(X_\gamma) $ 
  for 
  some $ \gamma<\beta $ thus $ P_e $ has already been defined, as well as the paths $ \{ P_b \}_{b\in G_\beta} $. Furthermore, the 
  path-system $ 
  \mathcal{P}_{\beta}:= \{P_e   \}_{e\in \mathsf{in}_D(X_\gamma)}\cup \{  P_b \}_{b\in G_\beta} $ is obviously edge-disjoint since any 
  two 
  elements of it are already members of $ \mathcal{P}_\gamma $ for some $ \gamma<\beta $.
   \end{sbiz}
  
  \begin{prop}\label{Mf két t-good uniója t-good}
  If $ X $ and $ Y $ are $ t $-good sets, then $ X\cup Y $ is  $ t $-good as well.
  \end{prop}
  
 \begin{sbiz}
  Let $ \mathcal{P}=\{ P_b \}_{b\in B_X}\cup \{ P_e \}_{e\in \mathsf{in}_D(X)} $ and $ \mathcal{Q}=\{ Q_b \}_{b\in B_Y}\cup \{ Q_e 
  \}_{e\in \mathsf{in}_D(Y)} $ be path-systems that show the $ t $-goodness of $ X $ and $ Y $ respectively.  Note that all the common 
  edges of the two path-systems are in $ A(D[X\cap Y]) $.  Let us define $ 
      B_Y'=\{ b\in B_{Y}: b\notin \mathsf{span}(B_X) \} $. For  \[ s\in B_Y'\cup \left[(\mathsf{in}_D(Y)\cap \mathsf{in}_D(X\cup 
      Y))\setminus 
      \mathsf{in}_D(X) \right],   \] let $ R_s $ be the path that we obtain by taking the initial segment of $ Q_s $ up to the first 
      vertex in $ X 
      $ 
      and concatenate it with $ P_e $ where $ e $ is the last edge of this terminal segment. The path-system    
  \[ \{ P_s:  s\in B_X\cup(\mathsf{in}_D(X)\cap \mathsf{in}_D(X\cup Y)) \}\cup \{ R_s:  s\in B_Y'\cup (\mathsf{in}_D(Y)\cap 
  \mathsf{in}_D(X\cup Y))\setminus 
          \mathsf{in}_D(X) \}     \]
          shows that $ X\cup Y $ is $ t $-good.
 \end{sbiz}
     
 \noindent Proposition \ref{Mf t-good lánc unió} and \ref{Mf két t-good uniója t-good} imply that the union of arbitrary many $ t 
 $-good sets is 
  $ t $-good thus the union of all of them (it is not an empty union because  $ \{ t \} $ is in it) as well.  
   \end{ssbiz}\\

   Our main tool to characterize the infeasible $ (i,e) $-extensions is the following theorem.
  \begin{thm}\label{Mf matroidos dualizált útfeladat feszítve}
 If the $ (I,t) $-linkage $ \mathcal{P}=\{ P_i \}_{i\in I} $ does not satisfy  the complementarity conditions with the largest $ t $-good 
 set $ T $, then there is a $ t 
 $-linkable $ I' $ for which $ \mathsf{span}(I') \supsetneq \mathsf{span}(I)$.
 
  \end{thm}

 \begin{biz}
 Assume that $ \mathcal{P} $ and $ T $ do not satisfy the complementarity conditions and for  $ i\in  I\setminus \mathcal{S}(T)=:I_{out} 
 $ the first edge of $ P_i $ that enters  $ T $ is $ e_i $. First we show that we may suppose without loss of generality 
  that there is a path-system $\{ P_i \}_{i\in B}\cup \{ P_e \}_{e\in \mathsf{in}_D(T)} $ such that
  
  \begin{enumerate}
  \item $\{ P_i \}_{i\in B}\cup \{ P_e \}_{e\in \mathsf{in}_D(T)} $ shows the $ t $-goodness of $ T $,
  \item $ B \subseteq I$ and $\{ P_i \}_{i\in B}\subseteq \mathcal{P} $,
  \item for $ i\in  I_{out} $ we have  $ A(P_i)\cap \mathsf{in}_D(T)=\{ e_i \} $,  and     $ 
   P_i[\mathsf{head}(e_i),t ]=P_{e_i} $. 
  \end{enumerate}
  
   Indeed, otherwise let $ J $ 
 be a maximal $ I/\mathcal{S}(T) $-independent subset  of $ 
  I_{out} $ and for $ j\in J $  take  
 the  segments $ \{ P_j[\mathsf{start}(P_j),\mathsf{head}(e_j)] \}_{j\in J}  $    from  $ \mathcal{P} $ and extend it to an $ (J \cup B, 
 t) $-linkage $ \mathcal{Q} $ by using the $ t 
  $-goodness of $ T $. Clearly $ I 
  \subseteq \mathsf{span}(J \cup B) $. We may assume that $ \mathsf{span}(I)=\mathsf{span}(J\cup B) $, otherwise $ I':=J\cup B $ would be 
  a 
  suitable choice for the theorem itself. We check 
  that $ A(\mathcal{Q})\cap\mathsf{in}_D(T)=\{ e_j \}_{j\in J}\subsetneq \mathsf{in}_D(T) $ by applying the fact that $ 
  \mathcal{P} $ and $ T $ do not satisfy 
  the complementarity conditions. Assume that the first complementarity condition fails for $ \mathcal{P} $ and $ T 
  $. We know $ \mathcal{S}(X)\subseteq \mathsf{span}(I) $  because of $ \mathsf{span}(I)=\mathsf{span}(J\cup B) $.  Thus  $ \{ e_j 
  \}_{j\in J} = \mathsf{in}_D(T)$ would mean that complementarity conditions  hold for $ \mathcal{P} $ and $ T $ which is not 
  the case. Finally the 
  edges $ \{ e_j 
    \}_{j\in J} \setminus \mathsf{in}_D(T)$  are unused by $ 
  \mathcal{Q} $, hence $ \mathcal{Q} $ and $ T $ do no satisfy the complementarity conditions either. If $ \mathcal{P} $ and $ T $ 
  satisfy the first complementarity condition, then by using the alternative formulation of complementarity conditions 
  2,3,4 (see at the end of Definition \ref{complementarity conditions})  we obtain that \[ \mathsf{in}_D(X)\setminus A_{last} (\{ 
  P_i[\mathsf{start}(P_i),\mathsf{head}(e_i)] \}_{i\in 
      I_{out}})\neq \varnothing.  \]  These edges will be unused by $ \mathcal{Q} $. 
  
    Now we turn to the proof of the theorem. Let us denote $ \{  i\in S: I+i\in  \mathcal{I} \}= I\cup (S\setminus \mathsf{span}(I))$ by 
    $ I^{\star} $.  We build an auxiliary digraph by extending $ D $. Pick 
    the  new vertices $ \{ u_i 
    \}_{i\in 
    I^{\star}},\  \{w_i \}_{i\in S} $ and $ s $ and draw the following 
additional 
edges

 \begin{enumerate}
\item $ \{ su_i: i\in I^{\star} \}, $
\item $ \{ u_iw_i: i\in I^{\star} \}, $ 
\item $ \{ w_iv: i\in S \wedge v\in \pi(i) \}. $
 
\end{enumerate}  
 
\noindent We denote the resulting digraph by $ \boldsymbol{D^{+}_0}=(V^{+},A_0) $. For $ i\in I $, we extend the path $ P_i $ with the 
new initial 
vertices $ 
 s,u_i, w_i $ to obtain the $ s 
\rightarrow t 
$ path $ P_i^{+} $ in $ D_0^{+} $. Let $ \boldsymbol{\mathcal{P}^{+}}=\{ 
P_i^{+}  \}_{i\in I} $. Finally, change the 
direction of the edges in $A(\mathcal{P}^{+})  $ to obtain $ \boldsymbol{D_0^{*}} $. We call these 
redirected edges 
\textbf{backward edges} and the others  \textbf{forward edges}. 
Let $ \boldsymbol{U_0^{+}} $ be the set of vertices of $ D_0^{*} $ that 
 are \textbf{u}nreachable  
 from $ s $ and let $ \boldsymbol{U_0}=U_0^{+}\cap V $. 
 
 Assume first that $ t\notin U_0 $ i.e. there is an $ s\rightarrow t $ path $ P^{+} $ in $ D_0^{+} $. Let its first edge be 
 $ su_{i_0} $. Note that $ i_0\in I^{\star}\setminus I $.  Use the standard 
 \textbf{augmentation path} technique to obtain a system of 
 edge-disjoit $ s\rightarrow t $ paths $ \{ Q_i^{+} \}_{i\in I+i_0} $ in $ D_0^{+} $ where the first edge of $ Q_i^{+} $ is $ su_i $. 
 More 
 precisely,  do the 
 following. First of all, keep 
 unchanged the paths $ P_i^{+}\in \mathcal{P}^{+} $  that have no 
 common edges with $ P^{+} $. After that,  take the symmetric 
 difference of $ A(P^{+}) $ and the united edge sets of the, 
 say $ k $ many, elements of $ 
 \mathcal{P}^{+} $ from which $ P^{+} $ uses some backward edges. From the 
 resulting edge set  build $ k+1 $ edge disjoint 
 $ s \rightarrow t $ paths  by the greedy method. Finally $ \{ Q_i^{+} \}_{i\in I+i_0} $ consists of the paths we kept unchanged and 
 these 
 $ 
 k+1 $ new paths.
 By cutting off the  three initial vertices   of the paths $ Q_i^{+} $,   we obtain a system 
 of edge-disjoint paths $\mathcal{Q}= \{ Q_i \}_{i\in I\cup \{ i_0 \}} $ 
 in $ D $ 
 such that for any  $ i\in I\cup \{ i_0 \} $, path $ Q_i $ goes from $ \pi(i) $ to $ t $, i.e. we get an $ (I+i_0,t) $-linkage. In this 
 case, 
 $ I':= I+i $ is appropriate.

Suppose that $ t\in U_0 $.  Clearly the paths in $ \mathcal{P}^{+} $ use all the edges in $ 
\mathsf{in}_{D_0^{+}}(U_0^{+}) $ and none of the edges in $ \mathsf{out}_{D_0^{+}}(U_0^{+}) $. Therefore the same holds for $ \mathcal{P} 
$ with 
respect to $ D $ and $ U_0 $. We claim  that $ T\subseteq U_0 $. Assume, to the contrary, that we have a strict $ s \rightarrow T $ path 
$ P^{+} $ 
in $ D_0^{*} $ 
with last edge $ f $. It follows from  
our additional assumptions about $ \mathcal{P} $ (the first paragraph of this proof) that it does not use any edge 
from $ \mathsf{out}_D(T) $, thus $ f $ cannot be a backward edge. If $ f\in \mathsf{in}_D(T)\setminus A(\mathcal{P}) $, then path $ 
P_f\in 
\{ P_e 
\}_{e\in \mathsf{in}_D(T)} $ would 
show the reachability of $ t $ from $ s $ in $ D_0^{*} $ contradicting  $ t\in U_0 $. Finally, suppose that  $f= w_iv $ for some $ v\in 
T $. Then $ i\in 
\mathcal{S}(T) $ and therefore $ i\in \mathsf{span}(B) $.  For 
$ j\in S\setminus I^{\star} $, the vertex  $ w_j $ has no ingoing edges in $ D_0^{*} $ hence we know that $ i\in I^{\star} $. Thus $ I+i  
$ 
is independent hence $ B+i $ as well. It follows that necessarily $ i\in B $.   But then the 
unique ingoing edge of $ w_i $ in $ D_0^{*} $ comes from  $ \mathsf{start}(P_i)\in T $   and it contradicts  the 
strictness 
of the $ s \rightarrow T $ path $ P^{+} $.  

Let $ \boldsymbol{I_{in,0}}=I\cap \mathcal{S}(U_0) $ and let $ I_{out,0}=I\setminus \mathcal{S}(U_0) $. We claim that for $ 
i\in I_{in,0} $ we have $\mathsf{start}(P_i) \in U_0 $. Indeed, otherwise  the backward edge $\mathsf{start}(P_i) w_i $ and 
any  forward edge $ w_iv $  with $ v\in \pi(i)\cap U_0 $ would lead to a contradiction with the definition of $ U_0^{+}  
$. 
Since $ A(\mathcal{P})\cap \mathsf{out}_D(U_0)=\varnothing $, we obtain that the path-system $ \{ 
P_i \}_{i\in I_{in,0}} $ lies in $ D[U_0] $. Then clearly the paths $ \{ P_i \}_{i\in I_{out,0}} $ have to use all the edges in $ 
\mathsf{in}_D(U_0) $ because $ \mathsf{in}_D(U_0)\subseteq A(\mathcal{P}) $. It follows that each of them uses exactly one such an edge, 
thus $ \mathcal{P} $ and $ U_0 $ satisfy all but possibly the first complementary conditions. 

Let us define    $ \boldsymbol{F_0} := \mathcal{S}(U_0)\setminus \mathsf{span}(I_{in,0}) $. Observe that $ F_0\neq\varnothing 
$, otherwise the first complementarity condition would hold for $ \mathcal{P} $ and $ U_0 $ and hence $ U_0 $ would be a $ t $-good set 
with $ 
U_0 \supsetneq T$ (clearly $  U_0\neq T $, since $ \mathcal{P} $ and $ T $ do not satisfy the complementarity conditions by assumption) 
which contradicts  the maximality of $ T $.

We know that $ \mathcal{S}(U_0)\subseteq \mathsf{span}(I) $, since for $ i\in 
S\setminus \mathsf{span}(I) $  the path  $ s, u_i, w_i,v $ where $ v\in \pi(i) $ shows that $ \pi(i)\cap U_0=\varnothing  $ i.e. $ 
i\notin \mathcal{S}(U_0) $. 
 Fix a well-ordering of $ I_{out,0} $. For $ i\in F_0 $, let $ s_0(i) $ be the smallest element of $ I_{out,0}\cap C(i,I) $.  Extend $ 
 D_0^{+} $ with the new edges $ \{ u_{s_0(i)}w_i: i\in F_0 \} 
$ to obtain 
$ D_{1}^{+}=(V^{+},A_1) $. We get $ 
D_1^{*} $ by changing the direction of edges in $ A(\mathcal{P}) $ in $ D_1^{+} $. Assume first that there is some $ s \rightarrow t $ 
path 
$ P^{+} $ in $ D_1^{*} $. For the first edge $ su_{i_0} $ of $ P^{+} $, we have $ i_0\in I^{\star}\setminus I $. Consider $ S_{rep}:=\{ 
i\in 
F_0: 
u_{s_0(i)}w_i\in A(P^{+}) \} $ and take the smallest element $ s_0(i) $ of $ s_0[S_{rep}] $. The set $ (I+i_0)-s_0(i)+i $ is independent, 
spans $ 
I+i_0 $, and the remaining elements of $ s_0[S_{rep}] $ have the same fundamental circuit on it as on $ I $. We can do recursively  in 
increasing order the other replacements  thus $I':= (I+i_0)+S_{rep}-s_0[S_{rep}] $ is independent and spans $ 
I+i_0 $. Applying  $ P^{+} $ in the augmentation path method  results in a desired $ (I',t) $-linkage. 

Assume that such a $ P^{+} $ does not exist.  Let 
$ U_1^{+}\ni t $ be the set of the vertices of $ D_1^{*} $ that are not reachable from $ s $ and let $ U_1=U_1^{+}\cap V $. Because of 
the 
new edges, 
$ U_1^{+} \subseteq U_0^{+} $ holds. Observe that the vertices  $ \{ w_i \}_{i\in \mathcal{S}(T)} $ did not get any new ingoing edge ($ B 
$ ensures $ \mathcal{S}(T)\cap 
F_0=\varnothing $) hence  $ T\subseteq U_1  $ follows in the same way as we proved $T \subseteq U_0 $. Let us define  $ I_{in,1}= I\cap 
\mathcal{S}(U_1) $ and $ I_{out,1}= I\setminus\mathcal{S}(U_1) $. The complementarity conditions 
hold for $ \mathcal{P} $ and $ U_1 $  except the first which may not,  and $ \mathcal{S}(U_1)\subseteq \mathsf{span}(I) $ holds. The 
proof of 
these 
facts are the same as for $ U_0 $.  Note that the new edges ensure that $ F_0\cap \mathcal{S}(U_1)=\varnothing $ hence for $ 
F_1:=\mathcal{S}(U_1)\setminus \mathsf{span}(I_{in,1}) $ we have $ F_0\cap F_1=\varnothing $. Let us 
extend 
the well-ordering of $ I_{out,0} $ to a well ordering of $ I_{out,1} 
$ in such 
a way that $ I_{out,1}\setminus I_{out,0} $ is a terminal segment in it. This choice ensures that for an edge $ u_{s_0(i)}w_i $ the 
element    $ 
s_0(i)$ 
is the 
smallest in $ 
I_{out,1}\cap C(i,I)   $, not just in $ I_{out,0}\cap C(i,I) $. For $ i\in F_1 $, let $ s_1(i) $ be the smallest element of $ 
I_{out,1}\cap C(i,I)   $. We obtain $ D_2^{+} $ from $ D_1^{+} $ by adding the new edges $ \{ 
u_{s_1(i)}w_i: i\in F_1  \} $.      

We define   the corresponding 
 notions $ D_2^{*},U_{2}^{+}, U_{2}, I_{in,2}, I_{out,2}, F_2 $ and continue the process recursively.  Suppose, to the contrary that we 
 do not 
 find a desired $ I' $.   Let us define $ D_\omega^{+}=(V^{+},\bigcup_{n<\omega}A_n) $ and the corresponding notions as earlier. Note 
 that $ 
U_\omega=\bigcap_{n<\omega} U_n \supseteq T $ and it satisfies all but  the first complementarity conditions (thus $ F_{\omega}\neq 
\varnothing $) with $ \mathcal{P} $.   Obviously  $ I_{in,\omega} \subseteq \bigcap_{n<\omega}I_{in,n} $ but in 
fact  $ I_{in,\omega} = \bigcap_{n<\omega}I_{in,n} $ holds. Indeed, $ i\in \bigcap_{n<\omega}I_{in,n} $ implies that $ P_i $ lies in $ 
\bigcap_{n<\omega} U_n=U_\omega $ and  $ \mathsf{start}(P_i) \in \pi(i) $  shows $ i\in I_{in,\omega} $. Let $ i\in F_\omega= 
\mathcal{S}(U_\omega)\setminus \mathsf{span}(I_{in,\omega}) $ be arbitrary. Then $ i\notin F_{n} $ for all $ n<\omega 
$ otherwise by the new edges we would have $ i\notin \mathcal{S}(U_{n+1}) \supseteq  \mathcal{S}(U_{\omega})$. On the other hand,   $ 
i\in 
\mathcal{S}(U_\omega) \subseteq \mathcal{S}(U_n) $ and by putting these together we obtain $ i\in \mathsf{span}(I_{in,n}) $ for all $ 
n<\omega $. We can not have $ i\in I_{in,n} $ for all $ n<\omega $ since then $ i\in I_{in,\omega} $ would follows. Suppose that $ 
i\notin I_{in,n} $ if $ n>n_0 $. Thus for any $n>n_0 $ we have $ (C(i,I)-i)\subseteq I_{in,n} $  but then  \[ (C(i,I)-i) \subseteq 
\bigcap_{n_0<n} I_{in,n}=\bigcap_{n<\omega} I_{in,n}= I_{in,\omega} \]  witnesses $ i\in \mathsf{span}(I_{in,\omega}) $,  hence $ i\notin 
F_{\omega} $ and thus $ 
F_{\omega}=\varnothing $ which is a 
contradiction.
\end{biz}\\

It worth  mentioning the following two consequences of the Theorem above.   
\begin{cor}
Let $ \mathcal{P} $ be a $ (I,t) $-linkage and assume that there exists a $ t $-linkable $ I' $ with $ \mathsf{span}(I') \supsetneq 
\mathsf{span}(I) $. Then one can choose such an $ I' $ and a $ (I',t) $-linkage $ \mathcal{Q} $ such that either $ 
A_{last}(\mathcal{Q})\subseteq A_{last}(\mathcal{P}) $ or the following hold:

\begin{enumerate}
\item $ \left|I'\setminus I\right|=\left|I\setminus I'\right|+1<\omega$,
\item $ A_{last}(\mathcal{P})\subseteq  A_{last}(\mathcal{Q})$,
\item $ \left|A_{last}(\mathcal{Q})\setminus 
 A_{last}(\mathcal{P})\right|=1 $,
 \item  $  \left|\mathcal{Q}\setminus \mathcal{P}\right|=\left|\mathcal{P}\setminus \mathcal{Q}\right|+1<\omega 
  $. 
\end{enumerate}
\end{cor}
\begin{sbiz}
Follow the proof of Theorem \ref{Mf matroidos dualizált útfeladat feszítve} without dealing with  $ T $  hence without the modification 
described in the first paragraph of that proof. Then we either obtain a desired $ (I',t) $-linkage that satisfies conditions 1-4 or a $ t 
$-good set $ 
U_n $  that satisfies the complementarity conditions with $ \mathcal{P} $. We may assume that the second 
possibility happens. Pick an arbitrary $ (J,t) $-linkage  $ \mathcal{R}=\{ R_i \}_{i\in J} $ with $ \mathsf{span}(J) \supsetneq 
\mathsf{span}(I) $. Let $ e_i $ be
the first 
edge of $ R_i $ that enters  $ U_n $, and let $ K $ be a maximal $ J/\mathcal{S}(U_n) $-independent subset of  $ J\setminus 
\mathcal{S}(U_n) $. Keep the segments $ R_i[\mathsf{start}(R_i),\mathsf{head}(e_i)] $ 
for $i\in K 
$. Concatenate $ R_i[\mathsf{start}(R_i),\mathsf{head}(e_i)] $ with $ P_j[\mathsf{head}(e_i),t] $, 
where $ 
P_j $ is the unique element of $ \mathcal{P} $ for which  $ A(P_j)\cap \mathsf{in}_D(U_n)=\{ e_i \} $, to obtain $ Q_i $, and let $ 
Q_i=P_i 
$  for $ i\in I\cap \mathcal{S}(X) 
$. It is routine to check that $ (I\cap \mathcal{S}(X))\cup K=:I' $ and $ \mathcal{Q}:=\{ Q_i  \}_{i\in I'} $ are 
appropriate.
\end{sbiz}
     
  \begin{cor}\label{Mf matroidos dualizált útfeladat}
   Suppose that $ I $ is $ t $-linkable.  Then there is no $ t $-linkable $ I' \supsetneq I $ if and only if
     there is a vertex set  $ X\ni t $ such that $ \mathcal{S}(X)\subseteq \mathsf{span}(I) $    
   and for $ I_{out}:=I \setminus\mathcal{S}(X) $ any strict $ (I_{out} ,X) $-linkage $ \mathcal{Q} $  we have $ 
   A_{last}(\mathcal{Q})=\mathsf{in}_D(X) $.
  \end{cor}
  \begin{sbiz}
  Apply Theorem \ref{Mf matroidos dualizált útfeladat feszítve} with the free matroid on $ I\cup (S \setminus \mathsf{span}(I)) $.
  \end{sbiz}\\

  Now we are able to prove the characterization of the infeasible $ (i,e) $-extensions. 
  \begin{lem}\label{Mf főlemma egyszerű}
  The $( i_0, e_0) $-extension is infeasible if and only if $ e_0 $ enters  some $ i_0 $-dangerous set.
  \end{lem} 
  
  \begin{sbiz}
  We have already checked the ``if'' so now we prove  the remaining direction.  Assume that vertex $ t $ witnesses the failure 
  of the linkage condition in the $ (i_0,e_0) $-extension 
  $ \mathfrak{R}_1 $ of $ \mathfrak{R} 
    $. 
   We claim that the largest $ t $-good set $ X $ with respect to $ \mathfrak{R}_1 $ is a desired $ i_0 $-dangerous set (with respect to 
   $ 
   \mathfrak{R} $). Let $ \mathcal{P}=\{ P_b 
   \}_{b\in B_{out}} $ be an arbitrary  reduced linkage for $ X $  with respect to $ \mathfrak{R} $, i.e. a strict $ (B_{out},X) 
   $-linkage   
    where $ 
   B_{out} $ is a base of $ \mathcal{N}(t)/\mathcal{S}(X) $. Note  that $ B_{out} $ contains a base $ B_{out}' $ of $ 
   \mathcal{N}(t)/\mathcal{S}_{\mathfrak{R}_1}(X) 
   $. Let $ \mathcal{P}'=\{ P_b 
      \}_{b\in B_{out}'} $.  Clearly $ e_0\in A(\mathcal{P}') $, otherwise from $ \mathcal{P}' $   one can get a $ (B,t) $-linkage with 
      respect to $ \mathfrak{R}_1 $ where 
   $ B $ is a base of $ 
   \mathcal{N}(t) $ by using the $ t $-goodness of $ X $. Suppose that $ e_0\in A(P_{b_1}) $. Then we are  able to construct a strict $ 
   (B-b_1,t) $-linkage $ \mathcal{L} $ with respect to $ \mathfrak{R}_1 $ from $ \mathcal{P}'\setminus \{ P_{b_1} \} $ via $ t $-goodness 
   as 
   above. By Fact 
   \ref{Mf bázis összehas}, an augmentation of this linkage in the sense of  
   Theorem 
     \ref{Mf matroidos dualizált útfeladat feszítve} would lead to a linkage for $ t  $ with respect to $ \mathfrak{R}_1 $, which is 
     impossible, therefore, by Theorem \ref{Mf matroidos dualizált útfeladat feszítve},   linkage $ \mathcal{L} $ satisfies the 
  complementarity 
  conditions with $ X $. Hence $ \mathcal{P}'\setminus \{ P_{b_1} \} $ needs to use all the edges in $ \mathsf{in}_{D-e_0}(X) $. The only 
  way 
  for  this to be true is if  $ e_0 $ is the last edge of $ P_{b_1},\  A_{last}(\mathcal{P}')=\mathsf{in}_D(X) $ and $ 
  \mathcal{P}'=\mathcal{P} $. Thus  $e_0\in \mathsf{in}_D(X)$ and $ X $ is tight with respect to $ \mathfrak{R} $. Furthermore, $ 
  \mathcal{P}'=\mathcal{P} $ implies $ \mathcal{S}(X)=\mathcal{S}_{\mathfrak{R}_1}(X) $, hence $ i_0\in \mathsf{span}(\mathcal{S}(X)) $. 
  Therefore, $ X $ is $ i_0 $-dangerous. 
  \end{sbiz}

\section{New matroid-rooted digraphs from tight sets}

 In finite combinatorics, it is a common proof technique to subdivide the problem into smaller sub-problems by using an appropriate 
 notion of tightness and then  solve the smaller sub-problems by induction from which one can obtain a solution for the original 
 problem. Unfortunately, in infinite combinatorics usually the resulting sub-problems are no longer ``smaller'' in any sense that would 
 make possible 
 such an induction. Even though do not lead to such an immediate success, the investigation of them could be fruitful, as happened in 
 this 
 topic.  

Through this chapter we have some fixed matroid-rooted digraph $ \mathfrak{R} $ that satisfies independence and the linkage condition.

\begin{claim}\label{Mf tightok struktúrája}
If $ X $ is a tight set  and $ \{ P_b \}_{b\in B_0} $ is a linkage for $ Z $  where  $\varnothing\neq  Z \subseteq X $, then 
\begin{enumerate}
\item the set $ B^{*}:=B_0\cap \mathcal{S}(X)$ is a base of $  \mathcal{N}(Z)\cap\mathsf{span}(\mathcal{S}(X)) $,
\item for $ b\in B^{*} $, path $ P_b $ lies in $ D[X] $,
\item for $ b\in B_0\setminus B^{*} $ we have $ \left|A(P_b)\cap \mathsf{in}_D(X)\right|=1 $,
\item all the edges $ \{ e\in \mathsf{in}_D(X): \mathsf{head}(e)\in \mathsf{to}_D(Z) \} $ are used by the paths $ \{ P_b \}_{b\in 
(B_0\setminus B^{*})} $.  
\end{enumerate}
\end{claim}

\begin{nbiz}
Pick a linkage $ \{ Q_b \}_{b\in B_1} $ for $ X $ and let $B_0 \subseteq B\subseteq B_0\cup B_1 $  a base of $ \mathcal{N}(X) $. Note 
that if $ A(P_{b_0})\cap A(Q_{b_1})\neq \varnothing $ for some $ b_0\in B_0 $ and $ b_1\in B_1 \setminus B_0 $, then $ b_1\in 
\mathcal{N}(Z)$ and 
therefore $ b_1\in \mathsf{span}(B_0) $. For $ b\in B\setminus B_0 $, let $ P_i=Q_i $, then  $\mathcal{P}':= \{ P_b \}_{b\in B} $ is a $ 
(B,X) 
$-linkage. 
 If $ B^{*}\subseteq B_0 $ is 
not a 
base of  
 \[ \mathcal{N}(Z)\cap\mathsf{span}(\mathcal{S}(X))\subseteq \mathsf{span}(B_0), \] then we may pick some $ i\in 
 (\mathcal{N}(Z)\cap\mathsf{span}(\mathcal{S}(X)))\setminus B^{*} $ for which $ B^{*}+i $ is independent and some  \[ j\in 
 C(i,B_0)\setminus B^{*}=C(i,B)\setminus B^{*}. 
 \]     
$ B-j+i $ is a base of $ \mathcal{N}(X) $ and  $ i\in \mathsf{span}(X) $ implies that for a suitable $ k\in \mathcal{S}(X) $ the set $ 
B-j+k $ as 
well (by Corollary 
\ref{Mf bázis-kör} with $ I:=B-j+i $ and $ J:=\mathcal{S}(X) $). Note that $ \mathsf{start}(P_{j}) \notin X $ because $ j\in B_0\setminus 
B^{*} $  and therefore $ A(Q_{i})\cap \mathsf{in}_D(X)\neq 
\varnothing $.  But then we may replace $ P_j $ by a trivial path $ P_k $ (consisting of a single vertex from $ \pi(k)\cap X $) in $ 
\mathcal{P}' $ 
and 
the new linkage does not use the edges 
$ A(P_{j})\cap \mathsf{in}_D(X)\neq \varnothing $ which contradicts  the tightness of $ X $.

If  for some $ b\in B^{*} $ path $ P_b $ is not entirely in $ X $, then it  
uses 
some element of $ \mathsf{in}_D(X) $. Replace $ P_b $ by a trivial path consisting of an element of $ \pi(b)\cap X $ to get a 
contradiction 
as above.

Assume that for some $ b\in B_0  $ path $ P_b $ uses more than one ingoing edge of $ X $, then we may replace it in $ \mathcal{Q} $ 
by 
its own initial segment up to 
the  head of its first edge  in $ \mathsf{in}_D(X) $ and get contradiction.

If   the linkage $\mathcal{P}= \{ P_b \}_{b\in 
B_0} $ does not use all the edges $ \{ e\in \mathsf{in}_D(X): \mathsf{head}(e)\in \mathsf{to}_D(Z) \} $, then the 
linkage $ \mathcal{P}'= \{ P_b \}_{b\in B} $ does not use these edges as well (since for $ b\in B\setminus B_{0} $ their heads may not 
even be reachable from $ \pi(b) $)  which contradicts  the tightness of $ X. \newmoon $
\end{nbiz}

\begin{cor}\label{Mf végtelen út nem elérhető}
Under Condition \ref{Mf végtelen előreutak} (Condition \ref{Mf végtelen hátrautak}), for a tight $ X $, a 
forward-infinite (backward-infinite) path $ P $ of $ D[X] $ may not be reachable in $ D $ from outside $ X $ (equivalently from $ \{ 
\mathsf{head}(e): e\in 
\mathsf{in}_D(X)  \} $).   
\end{cor}
\begin{sbiz}
Since  \[ \mathcal{N}(V(P))=\mathsf{span}(\mathcal{S}(V(P)))\subseteq \mathcal{N}(V(P))\cap \mathsf{span}(\mathcal{S}(X)), \]  by 
applying the first statement of  
Claim  \ref{Mf tightok struktúrája} 
with $ Z:=V(P) $ we obtain $ B^{*}=B_0 $ thus $ B_0 \setminus B^{*}=\varnothing $. Hence the Corollary follows from the fourth 
statement of Claim  \ref{Mf tightok struktúrája}.
\end{sbiz}\\
 
    For a tight $ X $, let  $ \boldsymbol{\mathfrak{R}[X]}$ be the matroid-rooted digraph with $ 
    D_{\mathfrak{R}[X]}=D[X],\   
   \mathcal{M}_{\mathfrak{R}[X]}= 
   \mathcal{S}(X)\bigoplus \{ i_e: e\in \mathsf{in}_D(X)\} $ where $ i_e $ are some new elements, distinct from the elements of 
   $ S $,  and we consider 
    $ \{ i_e: e\in \mathsf{in}_D(X)\} $ as a free matroid. Finally, let $ \pi_{\mathfrak{R}[X]}(i)=\pi(i)\cap X $ for $ 
   i\in 
   \mathcal{S}(X) $ and 
   let $ \pi_{\mathfrak{R}[X]}(i_e)=\{ \mathsf{head}(e)  \}$ for $ e\in \mathsf{in}_D(X) $. 
   
   \begin{obs}\label{Mf K, M matroidok kapcsolata}
    For $ U\cup\{ i \} \subseteq \mathcal{S}(X) $, we have   \[ i\in \mathsf{span}(\mathcal{S}(U))  \Longleftrightarrow i\in 
   \mathsf{span}_{\mathcal{M}_{\mathfrak{R}[X]}}(\mathcal{S}_{\mathfrak{R}[X]}(U)). \]   
   \end{obs}

\noindent Applying Claim \ref{Mf tightok struktúrája} we prove some basic facts related to $ \mathfrak{R}[X] $.

\begin{prop}\label{Mf dangerous teljesíti} 
\leavevmode

\begin{enumerate}
\item  $ \mathfrak{R}[X] $ satisfies the linkage condition and independence,
\item $ \mathsf{span}( \mathcal{N}(Z) \cap \mathcal{S}(X))=\mathcal{N}(Z)\cap \mathsf{span}(\mathcal{S}(X))\ (Z \subseteq X) $,
\item  
 $ \mathcal{N}_{\mathfrak{R}[X]}(Z)= \mathcal{N}(Z)\cap \mathcal{S}(X)\cup \{i_e: e\in 
\mathsf{in}_{D}(X)\wedge \mathsf{head}(e)\in \mathsf{to}_{D}(Z)   \}\ (Z\subseteq X)$.
 
\end{enumerate}  
\end{prop}
\begin{nbiz}
Let $ v\in X $ be arbitrary  and pick a linkage $ \{ P_b \}_{b\in B} $ for $ v $. Take the terminal segments of paths 
$ P_b $ from their first vertex in $ X $. Claim \ref{Mf tightok struktúrája} and the definition of $ \mathfrak{R}[X] $ ensure that the 
result is a linkage for $ v $ with 
respect to $ \mathfrak{R}[X] $. The independence preserving part follows from the fact that the  circuits of $ 
\mathcal{M}_{\mathfrak{R}[X]} $ are 
exactly those circuits of $ \mathcal{M} $ that lie in $ \mathcal{S}(X) $ and for $ Z\subseteq X $ we have $ \mathcal{S}(Z)=  
\mathcal{S}_{\mathfrak{R}[X]}(Z)\cap \mathcal{S}(X) 
$. Thus if we have a $ \mathcal{M}_{\mathfrak{R}[X]} $-circuit $ C \subseteq 
\mathcal{S}_{\mathfrak{R}[X]}(v) $ for some $ v \in X $, then $ C \subseteq \mathcal{S}(v) $ and $ C $ would be an $ \mathcal{M} 
$-circuit as 
well.

Assume that $ i\in \mathsf{span}( \mathcal{N}(Z) \cap \mathcal{S}(X)) $. Then by monotonicity   $ i\in \mathsf{span}(\mathcal{S}(X)) $ 
and 
$ i\in 
\mathsf{span}(\mathcal{N}(Z))=\mathcal{N}(Z) 
$, i.e. $ i\in  \mathcal{N}(Z)\cap\mathsf{span}(\mathcal{S}(X)) $. Suppose now $ i\in  \mathcal{N}(Z)\cap\mathsf{span}(\mathcal{S}(X)) $. 
By the first statement of Claim \ref{Mf tightok struktúrája}, we know that there is a base $ B^{*}\subseteq \mathcal{S}(X)\cap 
\mathcal{N}(Z) $ of $  \mathcal{N}(Z)\cap\mathsf{span}(\mathcal{S}(X)) $, hence

\[ i\in \mathsf{span}(B^{*})\subseteq \mathsf{span}( \mathcal{N}(Z) \cap \mathcal{S}(X)). \]

 At the third statement of this Proposition, the inclusion ``$ \subseteq $'' is straightforward. The linkage $ \{ P_b \}_{b\in B^{*}} $ 
 from the first two statements 
  of  Claim \ref{Mf tightok struktúrája} 
  ensures $ \mathcal{N}_{\mathfrak{R}[X]}(Z) \supseteq \mathcal{N}(Z)\cap \mathcal{S}(X) $ and the last two statements of Claim 
  \ref{Mf tightok struktúrája} show
 
 \[\mathcal{N}_{\mathfrak{R}[X]}(Z) \supseteq \{i_e: e\in 
 \mathsf{in}_{D}(X)\wedge \mathsf{head}(e)\in 
 \mathsf{to}_{D}(Z)  \}.\ \newmoon \]
  
\end{nbiz}

\begin{prop}
If $ X $ is a tight set and $ Z \subseteq X $, then $  Z $ is tight with respect to $ \mathfrak{R} $ if and only if  $ Z $ is tight with 
respect to 
$  \mathfrak{R}[X]$.
\end{prop}
\begin{sbiz}
Suppose that $ Z\neq \varnothing $ is not tight with respect to $  \mathfrak{R} $ and let $ \mathcal{P}=\{ P_b \}_{b\in B_0} $ be a 
linkage for $ Z $ 
such that for some $ f\in \mathsf{in}_D(Z) $ we have $ f\notin A(\mathcal{P}) $. Let $ q_b $ be the first vertex of $ P_b $ in $ X $.   
We show that the paths $ \{ 
P_b[q_b,\mathsf{end}(P_b)] \}_{b\in B_0} 
$ 
witnesses that $ Z $ is not tight with respect to $ 
\mathfrak{R}[X] $. By  Claim \ref{Mf tightok struktúrája} we know that these are really paths in $ D[X]  $.  Let $ B^{*}=\{ b\in B_0: 
\mathsf{start}(P_b)=q_b \} $. According to Claim \ref{Mf tightok struktúrája}   $ B^{*} $ is a base 
of $  \mathcal{N}(Z)\cap\mathsf{span}(\mathcal{S}(X)) $.  By the forth statement of Claim  \ref{Mf tightok struktúrája}, the set of the 
last edges of paths $ \{ 
 P_b[\mathsf{start}(P_b), q_b] \}_{ b\in B_0 
 \setminus 
B^{*}} $ is $ \{ e\in \mathsf{in}_D(X): \mathsf{head}(e)\in \mathsf{to}_D(Z) \}=:A_0 $. For  $ e\in A_0 $, let  $ 
P_{i_e}=  P_b[q_b,\mathsf{end}(P_b)] $ where  $ e $ is the last edge of 
 $ P_b[\mathsf{start}(P_b), q_b] $. The third statement of Proposition \ref{Mf 
dangerous teljesíti} ensures that the linkage

 \[ \{ P_b \}_{b\in B^{*}}\cup \{ P_{i_e} \}_{e\in A_0}  \]
 
\noindent corresponds to a base of $ 
\mathcal{N}_{\mathfrak{R}[X]}(Z) $. This linkage clearly does not use $ f\in \mathsf{in}_D(Z) $ but we are not done since we need to show 
$ f\in \mathsf{in}_{D[X]}(Z) $. Suppose, to the contrary that $ f\notin \mathsf{in}_{D[X]}(Z) $, then necessarily $ f\in 
\mathsf{in}_D(X)\cap \mathsf{in}_D(Z) $. But 
then by the last statement of Claim \ref{Mf tightok struktúrája} we obtain $ f\in A(\mathcal{P}) $ contradicting to the choice of $ f $. 
This completes  the proof of the ``only if'' part of the statement. 
 
 The proof of the other direction is very similar hence we give just a sketch. Take a linkage for $ Z $ which witnesses the untightness 
 of $ Z $ with 
 respect to $ \mathfrak{R}[X] $. Then give a 
 backward continuation for its paths in the form $ P_{i_e} $ by using an arbitrary linkage for $ Z $ with respect to $ \mathfrak{R} $. 
 The resulting linkage 
 for $ Z $ with respect to $ \mathfrak{R} $  shows the untightness of $ Z $ with respect to $ \mathfrak{R} $. 
\end{sbiz}\\

\noindent Observation \ref{Mf K, M matroidok kapcsolata} leads to the following consequence of the Proposition above. 
\begin{cor}\label{Mf dangerous kinn benn u.a.}
If $ X $ is a tight set,  $ Z \subseteq X $ and  $ i\in S $, then $ Z $ is $ i $-dangerous with respect to $ \mathfrak{R} 
$ if and only if $ Z $ is $ i $-dangerous with respect to $ \mathfrak{R}[X] $.
\end{cor}

\begin{claim}\label{Mf dangerousök metszete}
If $ X $ and $ Y $ are tight sets with $ X\cap Y\neq \varnothing $, then $ X\cap Y $ is tight as well. Furthermore if $ X $ and $ Y $ are 
$ i $-dangerous and $ i\in \mathcal{N}(X\cap Y) $, then $ X\cap Y $ is $ i $-dangerous.
\end{claim}

\begin{nbiz}
Let $ \{ P_b \}_{b\in B_0} $ be a linkage for $ X\cap Y=:Z $. If some edge enters  $ Z $, then it enters $ X $ or enters $ 
Y $ 
thus by applying the last statement of Claim \ref{Mf tightok struktúrája} to $ 
X $ with $ Z $ and then to $ Y $ with $ Z $ we obtain that the paths $ \{ P_b \}_{b\in B_0} $ use all the edges in $ 
\mathsf{in}_{D}(Z) $.

Let us turn to the dangerousness part of the Claim. By the first  statements of Claim \ref{Mf tightok struktúrája}, 
$ B^{*}_X:=B_0\cap \mathcal{S}(X) $ is a base of   $ \mathcal{N}(Z)\cap\mathsf{span}(\mathcal{S}(X)) $ and 
$ B^{*}_Y:=B_0\cap \mathcal{S}(Y) $ is a base of 
 $ \mathcal{N}(Z)\cap\mathsf{span}(\mathcal{S}(Y)) $. Clearly both $ B_X^{*} $ and $ B_Y^{*} $  need to contain a base of $ 
 \mathcal{S}(Z) $. This two bases of $ 
\mathcal{S}(Z) $ 
must be the same since $  B_X^{*}\cup  B_Y^{*} $ is independent. Therefore on the one hand  $  B_X^{*}\cap  B_Y^{*} $ contains a 
base of $ \mathcal{S}(Z) $. On the other hand, by the second statement of Claim \ref{Mf tightok struktúrája}  for $ b\in B_X^{*}\cap  
B_Y^{*} $ 
we have $ \mathsf{start}(P_b)\in X\cap Y=Z $ and hence $ b\in \mathcal{S}(Z) $  thus $  B_X^{*}\cap  B_Y^{*}\subseteq 
\mathcal{S}(Z) $. It follows that $ 
B_X^{*}\cap  B_Y^{*} $ is a base of $ \mathcal{S}(Z) $.  
  
 Assume now that $i\in \mathsf{span}( \mathcal{S}(X))\cap \mathsf{span}( \mathcal{S}(Y))  $ and $ i\in \mathcal{N}(Z) $ (hence $ i\in 
 \mathcal{N}(X)\cap \mathcal{N}(Y) $). Then  $ i\ \mathsf{span}( B^{*}_X)\cap  \mathsf{span}(B^{*}_Y) $. If $ i\in B^{*}_X\cap B^{*}_Y $, 
 then $ i\in \mathcal{S}(Z)$ and we are done. If exactly one element of $ \{ B^{*}_X, B^{*}_Y \} $ contains $ i $, then $ B^{*}_X\cup 
 B^{*}_Y $ 
 would contain a circuit through $ i $ which is impossible. Finally if 
 $ i $ 
 has a fundamental circuit on $ B^{*}_X $ and on $ B^{*}_Y $, then by  the independence of $ 
  B^{*}_X\cup B^{*}_Y $ and  the weak circuit elimination  (Fact \ref{weak circuit elimination})   
 these two circuits 
 must be the same and therefore lie in $ (B^{*}_X\cap B^{*}_Y)+i $. Hence   \[ i\in \mathsf{span}( B^{*}_X\cap B^{*}_Y)= 
 \mathsf{span}(\mathcal{S}(Z)).\ \newmoon \]   
\end{nbiz}

\begin{claim}\label{Mf belül megengedett kívül is}
If the $ (i,e) $-extension of $ \mathfrak{R}[X] $   is feasible where $ i\in S $, then the  $ (i,e) $-extension of $ \mathfrak{R} $ is 
feasible as well.
\end{claim}

\begin{sbiz}
Suppose, to the contrary, that it is not. Then by Lemma \ref{Mf főlemma egyszerű} $ e $ in an ingoing edge of some $ i $-dangerous set $ 
Y $. 
Using the fact that $ e $ lies in $ X $ we have 

 \[ i\in \mathcal{S}(\mathsf{tail}(e)) \subseteq \mathcal{S}(X)\subseteq \mathsf{span}(\mathcal{S}(X)) \]
 
\noindent hence $ X $ is $ i $-dangerous too. Edge $ e $ witnesses 
that $ i\in \mathcal{N}(X\cap Y) $, thus by Claim \ref{Mf dangerousök metszete}  $ Z:=X\cap Y $ is $ i $-dangerous with respect to $ 
\mathfrak{R} $ thus by Claim \ref{Mf dangerous kinn benn u.a.} it is $ i $-dangerous with respect to $ \mathfrak{R}[X] $ as well. But 
then the $ (i,e) $-extension of $ \mathfrak{R}[X] $ is infeasible since $ e\in \mathsf{in}_{D[X]}(Z)$ and $ Z $ is $ i $-dangerous which 
is 
a contradiction.
\end{sbiz}

\begin{cor}\label{külső bővítés meghatároz}
Let $ S' \subseteq S $. Then a feasible $ S' $-extension $ \mathfrak{R}[X]^{*} $ of $ \mathfrak{R}[X] $ of order $ n $ determines a 
unique, feasible  $ S' $-extension 
$ \mathfrak{R}^{*} $ of $ \mathfrak{R} $ of order $ n $ characterized by the property $ \mathfrak{R}^{*}[X]= \mathfrak{R}[X]^{*}$.
\end{cor}

\section{Augmentations at a prescribed vertex }
In this section we prove a Lemma, that allows us a kind of local augmentation. The Lemma will imply immediately  Theorem 
\ref{Mf főtétel} in the case of countable $ D $ and one can derive from it Theorem 
\ref{Mf főtétel} itself as well  without too much effort as we will do it in the last section. 

\begin{lem}\label{Mf egyet belerak könyebb}
 Assume that $ 
\mathfrak{R}=(D,\mathcal{M}, \pi) $ is independent and  satisfies the linkage condition and either Condition 
\ref{Mf végtelen előreutak} or Condition 
\ref{Mf végtelen hátrautak}. Then for any $ 
v\in V $ and for any $ W \subseteq S $ which is the union of finitely many components of $ \mathcal{M} $  there is a finite-order, 
feasible $ W $-extension $ 
\mathfrak{R}^{*} $ of $ 
\mathfrak{R} $  such that $ \mathcal{S}_{\mathfrak{R}^{*}}(v)\cap W $ is a base of $  \mathcal{N}(v)\cap W $. 
\end{lem}

\noindent  Assume, to the contrary, 
that the 
Lemma is false and choose an arbitrary  counterexample triple   $ 
\mathfrak{R}=(D,\mathcal{M},\pi),v_0, W $. 
We may assume (by replacing $ \mathfrak{R} $ by some feasible, finite-order $ W $-extension of itself) that we are not able to augmenting 
at 
$ v_0 
$ even by one. More precisely for any feasible, finite-order $ W $-extensions $ \mathfrak{R}' $ of $ \mathfrak{R} 
$ we have $ \mathcal{S}_{\mathfrak{R}'}(v_0)= \mathcal{S}_{\mathfrak{R}}(v_0)$.   
 Similarly we may suppose that $ \mathfrak{R} $ minimize the following expression among the feasible, finite-order  $ W $-extensions $ 
 \mathfrak{R}' $ of $ \mathfrak{R} $.
 
 \begin{align}\label{Mf legrövidebb út}
  \min \{ \left|A(P_{i_0})\right|: \{ P_i \}_{i\in B} \text{ is a reduced linkage for } v_0 \text{ with respect to }  
 \mathfrak{R}' 
 \text{  
 and } i_0\in B\cap W \} 
 \end{align}
 
Let the minimum for $ \mathfrak{R} $ be  taken on  $P_{i_0}\in \{ P_i \}_{i\in B} $. Consider the first edge $ e_0 $ of $ 
P_{i_0} $.

\begin{prop}\label{i,e nem megengedett útrövidülne}
The $ (i_0,e_0) $-extension of $ \mathfrak{R} $ is defined but not feasible.
\end{prop}
\begin{sbiz}
Suppose, to the 
contrary, that it is undefined i.e.  
$ i_0\in \mathsf{span}(\mathsf{head}(e_0)) $. Then by Corollary \ref{Mf bázis-kör} there is some $ i_0'\in 
\mathcal{S}(\mathsf{head}(e_0))\cap W\subseteq 
\mathcal{N}(v_0) $ such that $ B-i_0+i_0' $ is a base of $ \mathcal{N}(v_0)/\mathcal{S}(v_0) $ (which implies $ \mathsf{head}(e_0)\neq 
v_0 $ since $ \{ P_i \}_{i\in B} $ is a reduced linkage for $ v_0 $). But then we may  
replace $ i_0 $ by $ i_0' $ and $ 
P_{i_0} $ by  $P_{i_0'}:= P_{i_0}[\mathsf{head}(e_0), v_0] $  
 to get a contradiction with the fact that the minimum at (\ref{Mf legrövidebb út}) for $ \mathfrak{R} $ is $ \left|A(P_{i_0})\right| 
 $. On the other 
 hand,  the $ (i_0,e_0) $-extension cannot be feasible since otherwise the resulting extension would have a smaller minimum showed by the 
 linkage that we would obtain  from $ \{ P_i \}_{i\in B} $ by replacing $ P_{i_0} $ with $ P_{i_0}[\mathsf{head}(e_0), v_0] $.
\end{sbiz}\\ 

\noindent It follows by Lemma \ref{Mf főlemma egyszerű} that $ e_0 $ enters  some $ i_0 
$-dangerous set $ X $.

\begin{prop}\label{Mf v_0 nincs X-ben}
The set $ X $ does not contain $ v_0 $.
\end{prop}

\begin{sbiz}
Suppose, seeking for contradiction, that $ v_0\in X $. Then  all the paths in $ 
\{ P_{b} \}_{b\in B} $  meet $ X $.  
Pick a reduced linkage $ \{ P_b' \}_{b\in B_{X}} $ for $ X $  and note that 
if some 
$ P_b' $ have a common edge (or just a common vertex) with a path in $ \{ P_{b} \}_{b\in B} $, then $ b\in \mathsf{span}(B) $. 
Let $ B_{X}'=\{ b\in B_{X}: b\notin \mathsf{span}(B) \} $.   The set $ (B\cup B_{X}')\setminus\{ i_0 \} $ contains a base of $ 
\mathcal{N}(X)/\mathcal{S}(X) $ since $ B\cup B_{X}' $ clearly does and
$ i_0\in \mathsf{span}(\mathcal{S}(X)) $ implies that $ i_0 $ is a loop in  $ \mathcal{N}(X)/\mathcal{S}(X) $. The path-system 
$ (\{ P_{b} \}_{b\in B}\setminus \{ P_{i_0}\})\cup \{ P_b'\}_{b\in B_{X}'}  $ is edge-disjoint and shows that $ X $ is not tight, since 
the edge $ e_0\in \mathsf{in}_D(X) $ is unused, which is a contradiction. 
\end{sbiz}\\

\noindent By the first statement of Proposition 
\ref{Mf dangerous 
teljesíti}, $\mathfrak{R}[X] $ is independent and satisfies the linkage condition.

\begin{claim}
$ \mathfrak{R}[X] $ satisfies  Condition \ref{Mf végtelen előreutak} (Condition \ref{Mf végtelen hátrautak}).
\end{claim}  

\begin{nbiz}
Under  Condition \ref{Mf végtelen előreutak} the bases of $ \mathcal{M} $ are finite and hence a tight set may have just finitely many 
ingoing edges. Thus the bases of $ \mathcal{M}_{\mathfrak{R}[X]} $ are finite as well and therefore $ \mathfrak{R}[X] $ satisfies the 
first 
part 
of Condition \ref{Mf végtelen előreutak}.

In the case of Condition  \ref{Mf végtelen hátrautak},  observe that the independence of $ 
\mathfrak{R} $ implies that there are no loops in $ \mathcal{M} $. Thus from a component $ \mathcal{C} $ 
 of 
$ \mathcal{M} $ we get at most $ r(\mathcal{C}\cap \mathcal{S}(X))\leq r(\mathcal{C})<\infty $ 
 components of $ 
 \mathcal{M}_{\mathfrak{R}[X]} $.  
 Since under  Condition \ref{Mf végtelen hátrautak} the bases of $ \mathcal{M} $ are countable, a tight set may have just countably many 
ingoing edges. Hence the set of the further 
single-element components   $ \{ \{ i_e \}:  e\in \mathsf{in}_D(X)  \}$ is countable.  Thus $ \mathfrak{R}[X] $ satisfies the first part 
of Condition \ref{Mf végtelen hátrautak}.

 To show  the second part of Condition \ref{Mf végtelen 
előreutak} (Condition \ref{Mf végtelen 
hátrautak}) for $ \mathfrak{R}[X] $  take a forward-infinite (backward-infinite) path $ P $  of $ D[X] $ and let $ i\in 
\mathcal{N}_{\mathfrak{R}[X]}(V(P)) $ be arbitrary. By Corollary \ref{Mf végtelen út nem elérhető} and by the third statement of 
Proposition  \ref{Mf dangerous teljesíti} (with $ Z:=V(P) $), we obtain  \[i\in \mathcal{N}_{\mathfrak{R}[X]}(V(P)) 
=\mathcal{N}(V(P))\cap 
\mathcal{S}(X). \] 
Since Condition \ref{Mf végtelen 
előreutak} (Condition \ref{Mf végtelen 
hátrautak}) holds for $ \mathfrak{R} $ and $ i\in \mathcal{N}(V(P)) $ and Corollary \ref{Mf végtelen út nem elérhető} ensures   
\[ \mathcal{S}(V(P))=\mathcal{S}_{\mathfrak{R}[X]}(V(P))\subseteq \mathcal{S}(X),  \] we have  \[ i\in 
\mathsf{span}(\mathcal{S}(V(P)))=\mathsf{span}(\mathcal{S}_{\mathfrak{R}[X]}(V(P))). \]   
Hence 
by Observation \ref{Mf K, M matroidok kapcsolata} $ i\in 
\mathsf{span}_{\mathcal{M}_{\mathfrak{R}[X]}}(\mathcal{S}_{\mathfrak{R}[X]}(V(P))) $. \newmoon
\end{nbiz}\\
  
\noindent Let $ v_1 $ be the last vertex of $ P_{i_0} $ in $ X $. 

A restricted version of the following Sublemma was needed in 
\cite{fenyo2}  (applying the terminology of this paper $ \mathcal{M} $ was there the free matroid and $ \mathcal{N}(v)=\mathcal{M} $ for 
any $ v $). The more general circumstances make the precise formalisation of the proof a bit ugly although the new difficulties are just 
technical. After we finish the proof of the Sublemma (page \pageref{Sublemma vege}) we continue the proof of Lemma \ref{Mf egyet belerak 
könyebb}.

\begin{slem}\label{Mf Sublemma új ellenpélda}
$ \mathfrak{R}[X], v_1, W\cap\mathcal{S}(X)=:W^{*} $  is a counterexample for Lemma \ref{Mf egyet belerak könyebb}.
\end{slem}

\begin{sbiz}
Suppose, to the contrary, that it is not, and choose a feasible, finite-order $ W^{*} $-extension  $ \mathfrak{R}[X]^{*} $ of $ 
\mathfrak{R}[X] $ 
such that $ \mathcal{S}_{\mathfrak{R}[X]^{*}}(v_1)\cap W^{*}$  is a base of $\mathcal{N}_{\mathfrak{R}[X]}(v_1)\cap W^{*} $.  Then 
Corollary 
   \ref{külső  bővítés meghatároz}  gives a feasible, finite-order $ W^{*} $-extension $ \mathfrak{R}^{*} $ of $ \mathfrak{R} $ such that 
   $ 
    \mathfrak{R}^{*}[X]= \mathfrak{R}[X]^{*} $.   
 
 \begin{prop}\label{Mf i_0 generálódik v_1-ben könnyebb}
 There is some $ i_1\in \mathcal{S}_{ \mathfrak{R}^{*}}(v_1)\cap W^{*} $ for which 
 $B':= B-i_0+i_1 $ is a base of $ \mathcal{N}(v_0)/\mathcal{S}(v_0) $.
 \end{prop}
 \begin{sbiz}
  By the $ i_0 $-dangerousness of $ X $, we have  $ i_0\in \mathsf{span}(\mathcal{S}(X)) $ and path $ P_{i_0} $ shows $ i_0\in 
  \mathcal{N}(v_1) $.  Thus by applying Proposition \ref{Mf dangerous teljesíti} with $ Z:=\{ v_1 \} $ we 
 obtain
 
  \begin{align*}
  i_0 &\in\   \mathcal{N}(v_1)\cap \mathsf{span}(\mathcal{S}(X))\\
  &=\mathsf{span}(\mathcal{N}(v_1)\cap 
 \mathcal{S}(X))\\
 &=\mathsf{span}(\mathcal{N}_{\mathfrak{R}[X]}(v_1)\cap 
 \mathcal{S}(X)). 
  \end{align*} 
 Since $ W\ni i_0 $ and any circuit through $ i_0 $ lies in $ W $, it implies
  \[ i_0\in \mathsf{span}(\mathcal{N}_{\mathfrak{R}[X]}(v_1)\cap 
  \mathcal{S}(X)\cap W)=\mathsf{span}(\mathcal{N}_{\mathfrak{R}[X]}(v_1)\cap 
  W^{*})=\mathsf{span}(\mathcal{S}_{\mathfrak{R}[X]^{*}}(v_1)\cap W^{*}).  \]  
  
\noindent  Finally apply Corollary \ref{Mf bázis-kör} with $ 
  i:=i_0,\ I:=B, $ and $ J:=\mathcal{S}_{\mathfrak{R}^{*}}(v_1)\cap W^{*} $ and let $ i_1 $ be  the resulting $ j $ of Corollary \ref{Mf 
  bázis-kör} .
 \end{sbiz}

 \begin{claim}\label{Mf minimaliis hosszat csökkentő utrenszer}
 There is a base $ B^{*} $ of $ \mathcal{N}(v_0)/\mathcal{S}(v_0) $ and a  $ (B^{*},v_0) $-linkage $ \mathcal{P}^{*} $  with respect to $ 
\mathfrak{R}^{*} $ such that $ i_1\in B^{*} $ and $ 
   \mathcal{P}^{*}\ni P^{*}_{i_1}:=P_{i_0}[v_1,v_0] $. Hence (\ref{Mf legrövidebb út}) is smaller for $ \mathfrak{R}^{*} $ than for $ 
  \mathfrak{R} $ (which is a contradiction that proves Sublemma \ref{Mf Sublemma új ellenpélda}).
 \end{claim}
 
 \noindent Let $ B'= B-i_0+i_1 $ (see Proposition \ref{Mf i_0 generálódik v_1-ben könnyebb}).  If the paths $ \{ P_i \}_{i\in B'-i_1} $ 
 have no edge 
 in $A(D)\setminus A(D_{\mathfrak{R}^{*}})=:  \boldsymbol{A_{lost}} $, then $ \{ P^{*}_{i_1} \}\cup\{ P_i \}_{i\in B'-i_1} $ 
 shows that a desired linkage exists and we are done. We may assume that it is not the case.    Remember that $ A_{lost}\subseteq 
 A(D[X]) $.

Consider 
 the indices of those paths from $ \{ P_i \}_{i\in B'-i_1} $ that meet $ X $  i.e. $ 
 \boldsymbol{B_{ess}}:=\{i\in  B'-i_1: V(P_i)\cap X\neq \varnothing \} $ (the \textbf{ess}ential paths). Let us define the 
 \textbf{non}essential paths $ 
 B_{non}:=B'\setminus B_{ess} $ as well. For $ i\in B_{ess} 
 $ we denote by $ q_i 
 $ and $ z_i $ the first and the last vertex of $ P_i $ in $ X $ respectively. Whenever for some $ i\in B_{ess} $ the path $ 
 P_i[q_i,z_i] $ use an edge $ e\in \mathsf{out}_{D}(X) $ then there is an edge $ h_e\in \mathsf{in}_{D}(X) $ of $ P_i[q_i,z_i] $ which is 
 corresponding to the first ``come back'' to $ X $ after $ e $ (see Figure \ref{Mf H képe}). For all such an $ e $, we extend $ 
 D_{\mathfrak{R}^{*}}[X] $ with a new edge 
 $g(i,e) $ that goes from $ \mathsf{tail}(e) $ to $ \mathsf{head}(h_e) $. Furthermore pick a new vertex $ 
 t $ and for all  $ i\in B_{ess} $ draw an edge $ f_i $ 
from $ z_i $ to $ t $ to obtain $ \boldsymbol{H} $. 
 Let $ \boldsymbol{B_{in}}=\{ i\in B_{ess}: \mathsf{start}(P_i) \in X \} $ and let $ \boldsymbol{B_{out}}=B_{ess}\setminus B_{in} 
 $. 
 
 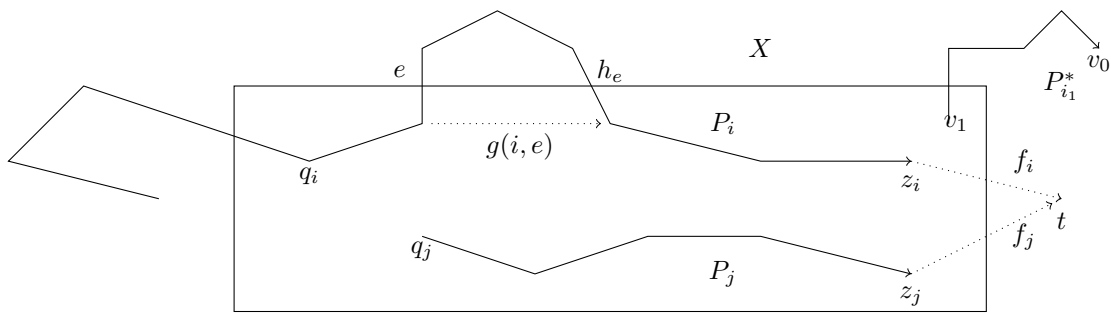
\begin{figure}[H]
 \centering
 
 \begin{tikzpicture}
 \draw (-4.5,3) rectangle (5.5,0);

 \draw[->] (-5.5,1.5) -- (-7.5,2) -- (-6.5,3) -- (-3.5,2) -- (-2,2.5) node (v2) {} -- (-2,3.5) -- (-1,4) -- (0,3.5) -- (0.5,2.5) node 
 (v3) {} 
 -- (2.5,2) -- (4.5,2);\draw[->,dotted] (4.5,2) -- (6.5,1.5) node (v1) {};
 
 \draw[->] (-2,1) -- (-0.5,0.5) -- (1,1) -- (2.5,1) -- (4.5,0.5);  \draw[->,dotted](4.5,0.5) -- (v1);
 
 \draw[->] (5,2.5) -- (5,3.5) -- (6,3.5) -- (6.5,4) -- (7,3.5);

 \node at (-3.5,1.8) {$q_i$};
 
 \node at (-2,0.8) {$q_j$};
 \node at (4.5,1.7) {$z_i$};
 
 \node at (4.5,0.2) {$z_j$};
 \node at (6,2) {$f_i$};
 \node at (6,1) {$f_j$};
 \node at (2.5,3.5) {$X$};
 
 \node at (6.5,3) {$P^*_{i_1}$};
 \node at (5.1,2.5) {$v_1$};
 \node at (7,3.3) {$v_0$};
 \draw[->,dotted]  (v2) edge (v3);
 \node at (-0.7,2.2) {$g(i,e)$};
 \node at (-2.3,3.2) {$e$};
 \node at (0.5,3.2) {$h_e$};
 \node at (6.5,1.2) {$t$};
 \node at (2,0.5) {$ P_j $};
 \node at (2,2.5) {$ P_i $};

 \end{tikzpicture}
 \caption{The construction of the digraph $ H $. We have $ i\in B_{out} $ and $ j\in B_{in} $.}\label{Mf H képe}
 \end{figure}

 We claim that one can justify Claim \ref{Mf minimaliis hosszat csökkentő utrenszer} by proving the following Claim.
 
 \begin{claim}\label{Mf belüljavít}
 There is  a $ B_{in}' \subseteq 
   \mathcal{S}(X)$ such that the set $( B' \setminus B_{in})\cup B_{in}'=B_{non}\cup B_{out}\cup B_{in}' $ is a base of $ 
   \mathcal{N}(v_0)/\mathcal{S}_{\mathfrak{R}}(v_0) $ and there is a system of 
   edge-disjoint 
   paths $ \{ 
  Q_i   \}_{i\in B_{out}\cup B_{in}'} $ in $ H $ such that for $ i\in B_{out} $ path $ 
  Q_i $ goes from $ 
  q_i $ to $ t $  and for $ i\in B_{in}' $ it goes from $ \pi_{\mathfrak{R}[X]^{*}}(i) $ to $ t $ .
 \end{claim}
 
\noindent    Indeed, for $ i\in 
 B_{non} $ let $ P^{*}_i=P_i $ if $ i\neq i_1 $ and let $ P^{*}_{i_1}=P_{i_0}[v_0,v_1] $. For $ i\in B_{out} $,  replace first the edges 
 in the form 
 $ g(j,e) $   of $ Q_i $ with 
 the corresponding path segments $ P_j[\mathsf{tail}(e),\mathsf{head}(h_e)] $. Then simplify the resulting walk to a path  and delete 
 its  last edge, say $ f_k $. Denote the result by  $ \widetilde{Q}_i $. Concatenate $ P_i[\mathsf{start}(P_i),q_i] $ with $ 
 \widetilde{Q}_i $ and   $ P_{k}[z_k,v_0] $ to obtain   
   $ P^{*}_i $. In the case $ i\in B_{in}' $,  we do the same, except we need to
 concatenate just $ \widetilde{Q}_i $ and  $ P_{k}[z_k,v_0] $ to get $ P^{*}_i $.  Finally  $  \{ P_i^{*} \}_{i\in B_{non}\cup 
 B_{out}\cup B_{in}'} $ 
 is  a desired linkage.

Let us define a matroid-rooted digraph that makes possible a reformulation of Claim \ref{Mf belüljavít}. For $ j\in B_{out} $, let  
$ F(j)=i_e \in S_{\mathcal{M}[X]} $ where $ e $ is the unique 
ingoing edge of $ q_j $ in $ P_j $ and let $ F $ be the identity on $ \mathcal{S}(X) $. We define 
$\boldsymbol{\mathcal{M}_{\mathfrak{Q}}}:=\left[\mathcal{S}(X)/(B_{non}\cup B_{out}\cup \mathcal{S}(v_0)) 
  \right]\bigoplus F[B_{out}] $.  Note that $ S_{\mathcal{M}_{\mathfrak{Q}}}\subseteq 
  S_{\mathcal{M}_{\mathfrak{R}[X]}} $.  Finally let $ 
  \boldsymbol{\mathfrak{Q}}:=(H,\mathcal{M}_{\mathfrak{Q}},\pi_{\mathfrak{R}[X]^{*}}\!|_{S_{\mathcal{M}_{\mathfrak{Q}}}} ) $.

\begin{obs}\label{Mf matroidok kapcsolata}
The $ \mathcal{M}_{\mathfrak{Q}} $-independent sets are $ \mathcal{M}_{\mathfrak{R}[X]} $-independent and  for $ S'\subseteq 
S_{\mathcal{M}_{\mathfrak{R}[X]}} $

\begin{equation}\label{Mf matroidok kapcsolata 1}
\mathsf{span}_{\mathcal{M}_{\mathfrak{R}[X]}}(S')\cap S_{\mathcal{M}_{\mathfrak{Q}}}\subseteq 
\mathsf{span}_{\mathcal{M}_{\mathfrak{Q}}}(S'\cap  S_{\mathcal{M}_{\mathfrak{Q}}} ).
\end{equation}  
 For any $ T\subseteq X $, we have   $ 
\mathcal{S}_{\mathfrak{Q}}(T) 
=\mathcal{S}_{\mathfrak{R}[X]^{*}}(T)\cap S_{\mathcal{M}_{\mathfrak{Q}}} $ which implies by using (\ref{Mf matroidok kapcsolata 1}) with 
$ 
S':= \mathcal{S}_{\mathfrak{R}[X]^{*}}(T) $

\begin{equation}\label{Mf matroidok kapcsolata 2}
\mathsf{span}_{\mathcal{M}_{\mathfrak{R}[X]}}(\mathcal{S}_{\mathfrak{R}[X]^{*}}(T))\cap 
S_{\mathcal{M}_{\mathfrak{Q}}}\subseteq 
\mathsf{span}_{\mathcal{M}_{\mathfrak{Q}}}(\mathcal{S}_{\mathfrak{Q}}(T)).
\end{equation}  
\end{obs}

\begin{prop}\label{Mf N Q és N R[X] }
For $ T\subseteq X $, we have $\mathcal{N}_{\mathfrak{R}[X]}(T)\cap S_{\mathcal{M}_{\mathfrak{Q}}}\subseteq \mathcal{N}_{\mathfrak{Q}}(T) 
$.
\end{prop}
\begin{nbiz}
We know that $ \mathcal{N}_{\mathfrak{R}[X]}=\mathcal{N}_{\mathfrak{R}[X]^{*}} $ since $ \mathfrak{R}[X]^{*} $ is a feasible extension of 
$ \mathfrak{R}[X] $. Obviously $ \mathsf{to}_{D_{\mathfrak{R}[X]^{*}}}(T)\subseteq \mathsf{to}_{H}(T) $ because $ 
D_{\mathfrak{R}[X]^{*}} $ is a subdigraph of $ H $. Then by applying (\ref{Mf matroidok kapcsolata 2}) of  Observation \ref{Mf matroidok 
kapcsolata} with $ T:=\mathsf{to}_{D_{\mathfrak{R}[X]^{*}}}(T) $
 \begin{align*}
  \mathcal{N}_{\mathfrak{R}[X]}(T)\cap 
  S_{\mathcal{M}_{\mathfrak{Q}}}&=\mathcal{N}_{\mathfrak{R}[X]^{*}}(T)\cap 
S_{\mathcal{M}_{\mathfrak{Q}}}=
\mathsf{span}_{\mathcal{M}_{\mathfrak{R}[X]^{*}}}\left[\mathcal{S}_{\mathfrak{R}[X]^{*}}(\mathsf{to}_{D_{\mathfrak{R}[X]^{*}}}(T)) 
\right] \cap
 S_{\mathcal{M}_{\mathfrak{Q}}}\\
 &\subseteq \mathsf{span}_{\mathcal{M}_{\mathfrak{Q}}} \left[\mathcal{S}_{\mathfrak{Q}}(\mathsf{to}_{D_{\mathfrak{R}[X]^{*}}}(T)) 
    \right]\subseteq\mathsf{span}_{\mathcal{M}_{\mathfrak{Q}}} \left[\mathcal{S}_{\mathfrak{Q}}(\mathsf{to}_{H}(T)) 
       \right]=\mathcal{N}_{\mathfrak{Q}}(T).\  \newmoon
 \end{align*}
 
\end{nbiz}

\noindent Clearly $\boldsymbol{B_0}:= F[B_{ess}]\subseteq \mathcal{N}_{\mathfrak{Q}}(t) $ is $ \mathcal{M}_{\mathfrak{Q}} $-independent. 
 In fact it is a base of $ \mathcal{N}_{\mathfrak{Q}}(t) $. Indeed, if there is an $ \mathcal{M}_{\mathfrak{Q}} $-independent $ I $  
 with   $ B_0\subsetneq I \subseteq 
 \mathcal{N}_{\mathfrak{Q}}(t) $, then we would obtain  \[ B' 
\subsetneq B_{non}\cup B_{out}\cup I\cap \mathcal{S}(X) \subseteq \mathcal{N}(v_0)/\mathcal{S}(v_0) \] where  $  B_{non}\cup B_{out}\cup 
I\cap \mathcal{S}(X) $ is $ \mathcal{N}(v_0)/\mathcal{S}(v_0) $-independent  which is impossible since $ B' $ is a base of $ 
\mathcal{N}(v_0)/\mathcal{S}(v_0) 
$.
 Thus an 
equivalent 
formulation of Claim \ref{Mf belüljavít}, that we will actually prove, 
 is the following. 

\begin{claim}\label{Mf nagyon belüljavít}
  There is a $ (\widehat{B},t) $-linkage with respect to $ \mathfrak{Q} $ where $ \widehat{B} $ is a base of  $ 
  \mathcal{N}_{\mathfrak{Q}}(t) $.
 
 \end{claim}
\begin{sbiz}
Fix a build sequence of $ \mathfrak{R}[X]^{*} $ from $ \mathfrak{R}[X] $ and let be the corresponding sequence of edges $ 
\left\langle h_m: m<M  \right\rangle $. Note that  $  \{ h_m: m<M\}= A_{lost} $. For $ n\leq M $, we denote the 
extension of $ H $ with the edges $ \{ h_m: n\leq m <M \} $ by $ H_n $ and we define $ 
\mathfrak{Q}_n=(H_n,\mathcal{M}_{\mathfrak{Q}},\pi_{\mathfrak{Q}}) $.   Note that $ H_M=H $ and hence $ \mathfrak{Q}_M=\mathfrak{Q} $.

\begin{obs}\label{Mf ob elérés úttal H D}
If $ H_n $ contains an $ u\rightarrow v $ path and $ v\neq t $, then $ D $ as well since we can just 
replace the edges in form $ g(j,e) $ 
by the corresponding paths of $ D $. 
\end{obs}
\begin{prop}
For all $ n\leq M $, we have $ \mathcal{N}_{\mathfrak{Q}_{n}}(t)=\mathcal{N}_{\mathfrak{Q}}(t) $.
\end{prop}

\begin{nbiz}
Obviously  $ \mathcal{N}_{\mathfrak{Q}_{n}}(t) 
\supseteq\mathcal{N}_{\mathfrak{Q}}(t)=\mathsf{span}_{\mathcal{M}_{\mathfrak{Q}}}(B_0) $. Suppose, to the contrary that $  
\mathcal{N}_{\mathfrak{Q}_{n}}(t) \setminus \mathcal{N}_{\mathfrak{Q}}(t)\neq \varnothing $. Then there is some $ i\in 
\mathcal{N}_{\mathfrak{Q}_{n}}(t)\setminus \mathsf{span}_{\mathcal{M}_{\mathfrak{Q}}}( B_0 )$ such that $ t $ (and hence $ \{ z_i 
\}_{i\in B_{ess}} $) is reachable from $ 
\pi_{\mathfrak{Q}}(i)  $ in $ H_n $. Necessarily $ i\in \mathcal{S}(X) $ because $ 
S_{\mathcal{M}_{\mathfrak{Q}}}\setminus \mathcal{S}(X)= F[B_{out}] \subseteq B_0 $. Then by Observation \ref{Mf ob elérés úttal H 
D}  $ 
\{ z_i \}_{i\in B_{ess}} $ is reachable 
from $ \pi_{\mathfrak{Q}}(i)=\pi_{\mathfrak{R}^{*}}(i) \cap X $  in $ D $. But then from $ \pi(i) $ as well, since all 
the new vertices that get  $ \pi(i) $ in an extension were originally reachable from $ \pi(i) $.
It 
follows that $ 
i\in \mathcal{N}(\{ z_i \}_{i\in B_{ess}})\subseteq \mathcal{N}(v_0) $  because $ v_0 $ is reachable in $ D $ from any 
element of $ \{ z_i \}_{i\in 
B_{ess}} $. But then $ B'\cup \{ i \}\subseteq \mathcal{N}(v_0)/\mathcal{S}(v_0) $ would be independent which is a contradiction since $ 
B' $ is a base of 
$ \mathcal{N}(v_0)/\mathcal{S}(v_0) $ and clearly $ i \notin B' $ since $ i\in S_{\mathcal{M}_{\mathfrak{Q}}} $ and \[ B'\cap 
S_{\mathcal{M}_{\mathfrak{Q}}}=B_{in}\subseteq B_0 
\subseteq
\mathsf{span}_{\mathcal{M}_{\mathfrak{Q}}}(B_0)\not\ni i.\ \newmoon \] 
\end{nbiz}

\noindent We 
prove by induction that for all $ n \leq M $ there is a $ (B_n,t) $-linkage with respect to $ \mathfrak{Q}_n $ where $ B_n $ is a base of 
$ 
\mathcal{N}_{\mathfrak{Q}}(t) $. For $ n=M $, we will obtain a desired linkage for Claim \ref{Mf nagyon belüljavít}. 

Let us start with the case $ n=0 $. For $ i\in B_{ess}  $, consider $ P_{i}[q_i,z_i] $ and for any $ e\in 
\mathsf{out}_D(X)\cap A(P_{i}[q_i,z_i])  $ replace the 
 segment $ P_i[\mathsf{tail}(e),\mathsf{head}(h_e)] $ by the single edge $ g(i,e) $  to obtain a path $ P_{F(i)}^{0} $ in $ H_0 $. The 
 linkage 
$\mathcal{P}_0= \{ 
P_{i}^{0} \}_{i\in B_0} $ is suitable.  

Suppose that there is a  $ (B_n,t) $-linkage $ \mathcal{P}_n=\{ P_i^{n} \}_{i\in B_n} $ with respect to $ 
\mathfrak{Q}_n $ such that $ n< M $ and $ B_n $ is a base of $ \mathcal{N}_{\mathfrak{Q}}(t) $. We need to 
give  a desired linkage with respect to   $ \mathfrak{Q}_{n+1} $. We may assume that for some $ j_0\in B_n $ we have $ 
h_n\in A(P_{j_0}^{n}) $ otherwise $ \mathcal{P}_n $ would be appropriate. Consider the $ (B_n-j_0,t) $-linkage $\mathcal{P}_{n}':= \{ 
P_i^{n} \}_{i\in 
B_n-j_0} $ in $ 
\mathfrak{Q}_{n+1} $. Note 
that if for some $ \mathcal{M}_{\mathfrak{Q}} $-independent $ I\subseteq 
\mathcal{N}_{\mathfrak{Q}_{n+1}}(t)=\mathcal{N}_{\mathfrak{Q}}(t) $ we have
\[  \mathsf{span}_{\mathcal{M}_{\mathfrak{Q}}}(I)\supsetneq 
\mathsf{span}_{\mathcal{M}_{\mathfrak{Q}}}(B_n-j_0), \] then by Fact \ref{Mf bázis összehas} $ I $ is necessarily a base of $ 
\mathcal{N}_{\mathfrak{Q}}(t) $. Apply Theorem \ref{Mf matroidos dualizált útfeladat feszítve} with 
linkage $ \mathcal{P}_{n}' $ in $ 
\mathfrak{Q}_{n+1} $.  We may assume that $ \mathcal{P}_{n}' $ and the largest $ t $-good set $ T^{+} $ of $ \mathfrak{Q}_{n+1} $ satisfy 
the 
complementarity conditions since otherwise Theorem \ref{Mf matroidos dualizált útfeladat feszítve} provides us a desired linkage. Let 
 $ f_{i(j_0)} $ be the last edge of $ P^{n}_{j_0} $.  Clearly  $ z_{i(j_0)}\in T^{+} $ otherwise $ 
f_{i(j_0)}\in \mathsf{in}_{H_{n+1}}(T^{+})\setminus A(\mathcal{P}_n') $ contradicting to the complementarity conditions. 

We build $ \mathcal{P}_{n+1} $ in three steps. First let $ \boldsymbol{B_n^{in}}=(B_n-j_{0})\cap \mathcal{S}_{\mathfrak{Q}}(T^{+}) 
$ and for $ 
i\in 
B_{n}^{in} $ let $ P^{n+1}_i=P_i^{n} $. By the first complementarity condition, $ B_n^{in} $ is a $ 
\mathcal{M}_{\mathfrak{Q}} $-base of $ 
\mathcal{S}_{\mathfrak{Q}}(T^{+}) $ and these paths lie inside $ T^{+} $. Let $ 
T=T^{+}-t $ and we define $ \boldsymbol{B_n^{un}}=B_n\setminus  \mathcal{N}_{\mathfrak{R}[X]}(T)$.

\begin{prop}
 $ j_{0}\notin B_n^{un} $.
\end{prop}

\begin{sbiz}
Since $ z_{i(j_0)}\in T $, the path $ P_{j_0}^{n} $ shows by applying 
Observation \ref{Mf ob elérés úttal H D}  that  $  T $  is 
reachable from $ \pi_{\mathfrak{Q}}(j_0)=\pi_{\mathfrak{R}[X]^{*}}(j_0) $ in $ D $ and hence from $ \pi_{\mathfrak{R}[X]}(j_0) $ as well 
thus 
$ j_0\in \mathcal{N}_{\mathfrak{R}[X]}(T) 
$.
 \end{sbiz}\\

\noindent In the second step we define $ P^{n+1}_i:=P^{n}_i $ for $ i\in B_n^{un} $.  Proposition above ensures that these paths are in $ 
H_{n+1} 
 $.   To construct the third 
part take a reduced linkage $ \mathcal{R}=\{ R_i \}_{i\in B_r} $ for $ T $ with respect to $ \mathfrak{R}[X]^{*} $. The path-system $ 
\mathcal{R} $ lies in $ H_{n+1} $ because  
 $ D_{\mathfrak{R}^{*}}[X] $ is a 
subdigraph of $ H_{n+1} $ . Since 
$ \mathfrak{R}[X]^{*} $ is a feasible extension of $ \mathfrak{R}[X] $, the set $ B_r $ is a base of $ 
\mathcal{N}_{\mathfrak{R}[X]}(T)/\mathcal{S}_{\mathfrak{R}[X]^{*}}(T) $.

\begin{prop} $  B_n \setminus 
(B_{n}^{in}\cup B_n^{un})\subseteq \mathsf{span}_{\mathcal{M}_{\mathfrak{Q}}}(B_n^{in}\cup (B_r\cap S_{\mathcal{M}_{\mathfrak{Q}}})) $. 
\end{prop}

\begin{sbiz}
 Let $ i\in B_n \setminus 
(B_{n}^{in}\cup B_n^{un}) $ be arbitrary. Then  $ i\notin \mathsf{span}_{\mathcal{M}_{\mathfrak{Q}}}(B_n^{in}) = 
\mathsf{span}_{\mathcal{M}_{\mathfrak{Q}}} 
(\mathcal{S}_{\mathfrak{Q}}(T)) $ hence by  (\ref{Mf matroidok kapcsolata 2}) of Observation \ref{Mf matroidok kapcsolata} $ i\notin 
\mathsf{span}_{\mathcal{M}_{\mathfrak{R}[X]}}(\mathcal{S}_{\mathfrak{R}[X]^{*}}(T)) 
$. On the other hand, $ i\in \mathcal{N}_{\mathfrak{R}[X]}(T) $ because $ i\notin B_n^{un} $. It shows that $ i\in 
\mathcal{N}_{\mathfrak{R}[X]}(T)/\mathcal{S}_{\mathfrak{R}[X]^{*}}(T) $. Then $ i\in\mathsf{span}_{\mathcal{M}_{\mathfrak{R}[X]}}(B_r\cup 
\mathcal{S}_{\mathfrak{R}[X]^{*}}(T)) $ because the choice of $ B_r $. 
Hence by  Observation \ref{Mf matroidok kapcsolata} $  i\in\mathsf{span}_{\mathcal{M}_{\mathfrak{Q}}}((B_r\cap 
S_{\mathcal{M}_{\mathfrak{Q}}})\cup \mathcal{S}_{\mathfrak{Q}}(T)) $ which is enough since $ 
\mathsf{span}_{\mathcal{M}_{\mathfrak{Q}}}(B_n^{in})=\mathcal{S}_{\mathfrak{Q}}(T) $. 
\end{sbiz}\\

\noindent We may take a $ B_r' \subseteq B_r\cap S_{\mathcal{M}_{\mathfrak{Q}}} $ for which $ B_n^{in}\cup B_n^{un}\cup B_r' $ is a 
maximal  $ 
\mathcal{M}_{\mathfrak{Q}} $-independent subset of $ B_n^{in}\cup B_n^{un}\cup (B_r\cap S_{\mathcal{M}_{\mathfrak{Q}}}) $.   For $ i\in 
B_r' $, concatenate $ R_i $ with the 
terminal segment of the  path $ P^{n}_j $ which is corresponding to the last edge of $ R_i $ to obtain $ P^{n+1}_{i} $. These paths also 
witnesses  that 
$ B_r'\subseteq \mathcal{N}_{\mathfrak{Q}}(t) $ and therefore $ B_n^{in}\cup B_n^{un}\cup B_r' $ is a base of $ 
\mathcal{N}_{\mathfrak{Q}}(t) $ since it is independent and spans such a base namely $ B_n $.

We need to check that the paths $ \{ P^{n+1}_{i} \}_{i\in B_r'} $ have no common edges with the paths  $ 
\{ P^{n}_i \}_{i\in 
B_n^{in}\cup 
B_n^{un}} $. The path-system   $ \{ P^{n}_{i} \}_{i\in B_n^{in}}  $ lies in $ T^{+} $ and the terminal segments of the paths $ \{ 
P^{n+1}_{i} 
\}_{i\in B_r'} $ from the first (and only)   entering to $ T^{+} $ are some other elements of $ \mathcal{P}_n $ which itself is an 
edge-disjoint 
system. Hence the path-system $ \{ P^{n+1}_{i} \}_{i\in B_r'}\cup \{ P^{n}_i \}_{i\in B_n^{in}}  $ is  edge-disjoint.  From the 
definition of $ B_n^{un} 
$  it  follows that

\[ \mathsf{to}_{D_{\mathfrak{R}[X]^{*}}}(T)\cap \bigcup_{i\in B_n^{un}}V(P^{n}_{i})=\varnothing. \]

On the other hand,

\[ \mathsf{to}_{D_{\mathfrak{R}[X]^{*}}}(T) \supseteq \bigcup_{i\in B_r}V(P_i^{n+1})\setminus \{ t \}, \]

thus the two paths-systems may not even have a common vertex other than $ t $. 
Now the proof of Claim \ref{Mf nagyon belüljavít} is 
complete and hence the proof of Claim \ref{Mf minimaliis hosszat csökkentő utrenszer} and the proof of Sublemma \ref{Mf Sublemma új 
ellenpélda} as well.\label{Sublemma vege}  
\end{sbiz}
\end{sbiz}\\

We continue the proof of Lemma \ref{Mf egyet belerak könyebb}. We obtained by Sublemma \ref{Mf Sublemma új ellenpélda} and by Proposition 
\ref{Mf v_0 nincs X-ben}  that if  $ \mathfrak{R}_0, v_0, W_0 $  is a counterexample triple, 
then there is a feasible, finite-order $ W_0 $-extension $ 
\mathfrak{R}_1 $ of $ \mathfrak{R}_0 $ such that there is a  vertex set $ X=:X_1\not\ni v_0 $ which is tight wit respect to $ 
\mathfrak{R}_1 $ and 
for a 
suitable 
$v_1\in  X_1 $ the triple $ \mathfrak{R}_1[X],v_1, \mathcal{S}_{\mathfrak{R}_1}(X)\cap W_0=:W_1  $ is a counterexample again. Furthermore 
we 
know that 
there is an $ e_1\in \mathsf{in}_{D_{\mathfrak{R}_1}}(X_1) $ and there is a path, namely $ P_{i_1}^{*} $ see Figure \ref{Mf H képe}, that 
goes 
strictly 
from $ X_1 $ to $ v_0 $ and starts at $ v_1 $. The path $ P_{i_0} $ shows that $ v_1 $ is reachable outside $ X_1 $ in $ 
D_{\mathfrak{R}_1} 
$.  We may apply these observations with the new counterexample triple and iterate the process 
recursively to get an infinite sequence of counterexample triples $ \left\langle (\mathfrak{R}_n[X_n], W_n, v_n): n<\omega   
\right\rangle$ with $ X_0:=V $. 

Here $ \mathfrak{R}_{n+1} $ is a finite-order, feasible $ W_n $-extension of 
$ \mathfrak{R}_n $ where the extension use edges only from $ D[X_n] $, $ \left\langle X_n: n<\omega  \right\rangle $ is a nested sequence 
of 
vertex sets such that $ X_n $ is tight with respect to $ \mathfrak{R}_n $ and $ v_n\in X_n $ but $ v_n\notin X_{n+1} $. We also have a 
path $ P_n  $ in $ D_{\mathfrak{R}_n} $ from $ X_{n+1} $ to $ v_{n} $ with $  \mathsf{start}(P_n)=v_{n+1} $ and some edge $ e_{n+1}\in 
\mathsf{in}_{D_{\mathfrak{R}_{n+1}}[X_n]}(X_{n+1}) $.

If $ \mathfrak{R}_0 $ satisfies Condition \ref{Mf végtelen hátrautak}, then we build a backward-infinite path $ P $ by concatenating the 
paths $ P_n $ for  $ n=1,2,\dots $. Then $ P $ lies in the $ \mathfrak{R}_1 $-tight $ X_1 $ and $ V(P) $  is reachable in $ 
D_{\mathfrak{R}_1} $ from outside $ X_1 $ in $ D_{\mathfrak{R}_1} $ (showed by $ P_{i_1} $) contradicting to Corollary \ref{Mf végtelen 
út nem elérhető}. It proves Lemma \ref{Mf egyet belerak könyebb} in the case when $ \mathfrak{R}_0 $ satisfies Condition \ref{Mf végtelen 
hátrautak}.

Suppose now that $ \mathfrak{R}_0 $ satisfies Condition \ref{Mf végtelen előreutak}. The sequence $ \left\langle \mathcal{N}(X_n): 
n<\omega  \right\rangle $ is  $ \subseteq $-decreasing  thus  $ \left\langle r(\mathcal{N}(X_n)): 
n<\omega  \right\rangle $ is a decreasing sequence of natural numbers therefore  by throwing away some initial elements we may assume 
that $ r(\mathcal{N}(X_n)) $ does not depend on $ n $. On the other hand,  the $ 
(n+1) $-th 
extension uses edges only from $ D_{\mathfrak{R}_n}[X_n] $ thus we have 

\[ \mathcal{S}_{\mathfrak{R}_n}(X_n)=\mathcal{S}_{\mathfrak{R}_{n+1}}(X_n)\supseteq \mathcal{S}_{\mathfrak{R}_{n+1}}(X_{n+1}) \] and 
therefore \[ r(\mathcal{S}_{\mathfrak{R}_n}(X_n)) \geq r(\mathcal{S}_{\mathfrak{R}_{n+1}}(X_{n+1})). \] But 
then   
\[ r(\mathcal{N}(X_n)/\mathcal{S}_{\mathfrak{R}_n}(X_n))=r(\mathcal{N}(X_n))-r(\mathcal{S}_{\mathfrak{R}_n}(X_n)) \] is an increasing 
function of $ n $  (bounded by $ 
r(\mathcal{M})<\infty $) hence similarly  we may suppose that it is constant, say $ m_0 $. Pick a reduced linkage $ \mathcal{P}_1 $ for $ 
X_{1} $ with respect to $ 
\mathfrak{R}_{1} $. It consists of $ m_0 $ paths and these paths use all the elements of $ \mathsf{in}_{D_{\mathfrak{R}_1}}(X_1)\ni e_1 $ 
because of the tightness of $ X_1 $. Then pick a reduced linkage $ \mathcal{Q} $ for $ X_{2} $ in $ 
\mathfrak{R}_{2} $. Observe that these paths also use all the elements of $ 
\mathsf{in}_{D_{\mathfrak{R}_1}}(X_1)=\mathsf{in}_{D_{\mathfrak{R}_2}}(X_1) $. Take the set of the terminal segments of the elements of $ 
\mathcal{Q} $  
from the first vertex in $ X_1 $ and denote it by $ \mathcal{Q}' $.  From $ \mathcal{P}_1 $ obtain via concatenation with elements in $ 
\mathcal{Q}' $ a reduced linkage $ \mathcal{P}_2 $ for $ X_2 $ with respect to $ \mathfrak{R}_2 $.  Iterate the process recursively. In a 
general step we have a reduced 
linkage $ \mathcal{P}_n $ for $ X_{n} $ with respect to $ \mathfrak{R}_n $ and we find forward-continuations for the elements of $ 
\mathcal{P}_{n} $ to obtain a reduced linkage for $ X_{n+1} $ with respect to $ \mathfrak{R}_{n+1} $. By the tightness of $ X_{n+1} $ 
with respect to $ \mathfrak{R}_{n+1} $,  necessarily  $ e_{n+1}\in A(\mathcal{P}_{n+1}) $. Eventually we obtain  an edge disjoint 
path-system $ 
\mathcal{P} $ with $ 
m_0 $ members. Since the edges $ \{ e_n \}_{1 \leq n <\omega}\subseteq A(\mathcal{P}) $ are pairwise distinct, there is a $ P\in 
\mathcal{P} $ that contains infinitely many of them. A terminal segment of the forward-infinite path $ P $ lies inside $ X_1 $  and 
reachable from outside $ X_1 $ in $ D_{\mathfrak{R}_1} $ (showed by $ P $ itself) which contradicts  Corollary 
\ref{Mf végtelen út nem elérhető}. Now the proof of Lemma \ref{Mf egyet belerak könyebb} is complete. \rule{1.5ex}{1.5ex}

\section{Careful iteration of local augmentations}

Now we are able to prove our main result Theorem \ref{Mf főtétel}. Suppose first that $ V $ is countable and  $ V=\{ v_n \}_{n<\omega} $ 
and the components of $ \mathcal{M} $ are $ \{ \mathcal{C}_n \}_{n<\omega} $ (if there are just finitely many, then repetition is 
allowed). 
Let $ \omega \times \omega= \{ p_n: n<\omega \} $. We build recursively a sequence $ \left\langle \mathfrak{R}_n: n\leq\omega  
\right\rangle $ such 
that 
$ \mathfrak{R}_0=\mathfrak{R} $ and if $ p_n=\left\langle m,k  \right\rangle $, then we obtain $ \mathfrak{R}_{n+1} $ by applying Lemma 
\ref{Mf egyet belerak könyebb} to $ \mathfrak{R}_n=(D_n, \mathcal{M}, \pi_n) $ with $ v_m $ and $ \mathcal{C}_k $. Finally let $ 
\mathfrak{R}_\omega= 
(D_\omega, \mathcal{M},\pi_\omega) $ where $ D_\omega=(V, \bigcap_{n<\omega} A(D_n)) $ and $ \pi_{\omega}(i)=\bigcup_{n<\omega}\pi_n(i) 
$. By the construction, for any $ v\in V $ and any component $ \mathcal{C} $ of $ \mathcal{M} $,  the set $ \mathcal{S}_\omega(v)\cap  
\mathcal{C}
$ is a base of $ \mathcal{N}(v)\cap  \mathcal{C} $ thus for all $ v\in V $ the set $ \mathcal{S}_\omega(v) $ is a base of $ 
\mathcal{N}(v) $.

In the general case, we should organize the recursion more warily to ensure that after limit steps Condition \ref{Mf fenntartandó 
csúcsra} 
holds. Let $ 
V=\{ v_\xi: \xi<\kappa \} $. To obtain $ \mathfrak{R}_{\xi+1} $ from $ \mathfrak{R}_{\xi} $  we consider $ v_\xi $ and all the finitely 
many vertices that lost some ingoing edge since the last limit step. We apply to these vertices $ v $ one by one 
in an arbitrary 
order  Lemma \ref{Mf egyet belerak 
könyebb} with 
the smallest $ n $ for which $ \mathcal{S}_\xi(v)\cap  \mathcal{C}_n $ is not a base of $ \mathcal{N}(v)\cap \mathcal{C}_n $ (if such an 
$ n $ does not exists for some $ v $, then do nothing with that $ v $). Observe 
that it ensures that after a limit step $ \alpha $ a $ v\in V $ either keeps all of its ingoing edges or $ \mathcal{S}_\alpha(v) $ is a 
base of $ \mathcal{N}(v) $. We need to justify that at limit steps we obtain feasible extensions in the process above.

\begin{prop}
Let $ \alpha<\kappa $ be a limit ordinal and suppose that $ \left\langle \mathfrak{R}_\beta: \beta<\alpha  \right\rangle $ has been 
 defined as above and this is a chain of feasible extensions of $ \mathfrak{R}_0=\mathfrak{R} $. Then the 
 limit  $ 
\mathfrak{R}_\alpha $ 
of the sequence is also feasible extension of $ \mathfrak{R} $.
\end{prop}
\begin{sbiz}
Let $ v\in V $ arbitrary and pick a linkage $ \{ P_b \}_{b\in B} $ for $ v $ with respect to $ \mathfrak{R} $. If some $ P_b $ is not 
a path in $ D_\alpha $ then replace it by the terminal segment $ Q_b $ of itself  that starts at the head $ u_b $ of the last deleted 
edge 
of $ P_b $ otherwise let $ Q_b=P_b $ and let $ u_b  $ be the first vertex of $ P_b $. Note that the our recursive process guarantees that 
$ b\in 
\mathsf{span}(\mathcal{S}_{\alpha}(u_b)) 
$. It is enough to show that there is a transversal for $ 
\{ \mathcal{S}_\alpha(u_b) \}_{b\in B} $ which is a base of $ \mathcal{N}(v) $. To do so we prove that for any component $ \mathcal{C} $ 
of $ \mathcal{M} $ there is a 
transversal for $ \{ \mathcal{S}_\alpha(u_b) \}_{b\in B\cap  \mathcal{C}} $ which is a base of $ \mathcal{N}(v)\cap  \mathcal{C} $. 

Let $ \mathcal{C} $ be fixed and let $ 
B_{\mathcal{C}}:=B\cap \mathcal{C}=\{ b_1,\dots ,b_{\ell_0} \} $. Pick a base $ B_{\mathcal{C}}'=\{ b_1',\dots,b_{\ell_0}' \} $ of $ 
\mathcal{N}(v)\cap  \mathcal{C} $ for which 
$ b_\ell'\in 
\mathsf{span}(\mathsf{S}_{\alpha}(u_{b_\ell}))  $ holds for all $ 1\leq \ell \leq \ell_0 $ and  $ b_\ell'\in 
\mathcal{S}_{\alpha}(u_{b_\ell})  $ for as many $ \ell $  as possible. Assume, to the contrary, that $ 
b_{\ell_1}'\notin 
\mathcal{S}_{\alpha}(u_{b_{\ell_1}})  $ for some $ 1 \leq \ell_1 \leq \ell_0 $.  The fact $ b_{\ell_1}'\in 
\mathsf{span}(\mathcal{S}_{\alpha}(u_{b_{\ell_1}}))\setminus \mathcal{S}_{\alpha}(u_{b_{\ell_1}}) $ implies that there is a circuit $ 
C\ni 
b_{\ell_1}' $ such that $ 
(C\setminus \{ b_{\ell_1}' \}) \subseteq \mathcal{S}_{\alpha}(u_{b_{\ell_1}})  $. Note that $ C \subseteq  \mathcal{C} $ because $ 
b_{\ell_1}'\in 
 \mathcal{C} $. Since $(C\setminus \{ b_{\ell_1}' \})\not\subseteq \mathsf{span} (B_{\mathcal{C}}'-b_{\ell_1}') $ (otherwise $ 
 b_{\ell_1}'\in 
 \mathsf{span} 
(B_{\mathcal{C}}'-b_{\ell_1}') $ ), there is some $ b_{\ell_1}''\in (C\setminus \{ b_{\ell_1}' \}) $ for which $ B_{\mathcal{C}}'- 
b_{\ell_1}'+ b_{\ell_1}'' $ is still a base of $ \mathcal{N}(v)\cap  \mathcal{C} $ contradicting to the choice of $ B_{\mathcal{C}}' 
$.
\end{sbiz}

\end{document}